\documentclass{amsart}
\usepackage{graphicx}
\newtheorem{thm}{Theorem}[section]

\newtheorem{lem}[thm]{Lemma}
\newtheorem{prop}[thm]{Proposition}
\theoremstyle{definition}

\theoremstyle{remark} \theoremstyle{Proof}
\newtheorem{rem}[thm]{Remark}

\numberwithin{equation}{section}

\author[O. Ajebbar]{ Ajebbar Omar}
\address{ Ajebbar Omar\\Department of Mathematics and Computer Sciences\\
Sultan Moulay Slimane University, Multidisciplinary faculty, Beni Mellal\\
Morocco} \email{omar-ajb@hotmail.com}
\author[E. Elqorachi]{ Elqorachi  Elhoucien}
\address{ Elqorachi Elhoucien\\Department of Mathematics\\
Ibn Zohr University, Faculty of Sciences, Agadir\\
Morocco} \email{elqorachi@hotmail.com}

\title[A system of cosine-sine functional equations]{A system of cosine-sine functional equations on a semigroup generated by its squares}

\begin{document}
\begin{abstract}
Given a semigroup $S$ generated by its squares,
we determine the complex-valued solutions of the following system of cosine-sine functional equations
\begin{align*}
f(xy)=f(x)g_{1}(y)+g_{1}(x)f(y)+\lambda_{1}^{2}\,h(x)h(y),\; x,y\in S,\\
h(xy)=h(x)g_{2}(y)+g_{2}(x)h(y)+\lambda_{2}^{2}\,f(x)f(y),\; x,y\in S,
\end{align*}
where $\lambda_{1},\lambda_{2}\in\mathbb{C}$ are given constants and $f,g_{1},g_{2},h:S\to\mathbb{C}$ are unknown functions.
\end{abstract}
\maketitle
\section{Introduction}
Let $S$ be a semigroup. The cosine-sine functional equation is
\begin{equation}\label{Eq8-0-0}f(xy)=f(x)g(y)+g(x)f(y)+h(x)h(y),\;x,y\in S,\end{equation} where $f,g,h:S\to\mathbb{C}$. This equation was solved by Chung, Kannappan and Ng \cite{C-K-N} for the case that $S$ is a group.
The result in \cite{C-K-N} were extended by the authors \cite{Ajb and Elq1} to the case that $S$ is a semigroup generated by its squares.
The functional equation (\ref{Eq8-0-0}) generalizes  the sine addition formula
\begin{equation}\label{Eq8-0-1}f(xy)=f(x)g(y)+g(x)f(y)\;x,y\in S,\end{equation}
and the cosine addition formula
\begin{equation}\label{Eq8-0-2}f(xy)=f(x)f(y)+g(x)g(y)\;x,y\in S.\end{equation}
The most current result on semigroups was given recently by Ebanks \cite{Ebanks1}, \cite{Ebanks} and \cite{Ebanks2}. Recently, Ebanks \cite{COSINE} find the solutions of the functional equation (\ref{Eq8-0-0}) on a large class of semigroups and discuss the obstacles to finding a general solution for all semigroups.\\
The functional equations (\ref{Eq8-0-1}), (\ref{Eq8-0-2}) has been studied extensively (both individually and as part of a system of equations). See for example \cite{Acz}, \cite{Acz-Dh} and \cite{Stetk1}.\\ The approach in \cite{C-K-N} is based on solving a system of three functional equations (See pages 268-271). So, it is quite naturel that a more general functional equation than (\ref{Eq8-0-0}) requires more computation and leads to a Levi-Civita type system of functional equations. The linear independence of the solutions is an important phenomenon for obtaining such systems.\\
In this paper we dealt with the following system of Levi-Civita functional equations
\begin{equation}\label{eq01}f(xy)=f(x)g_{1}(y)+g_{1}(x)f(y)+\lambda_{1}^{2}\,h(x)h(y),\;x,y\in S,\end{equation}
\begin{equation}\label{eq02}h(xy)=h(x)g_{2}(y)+g_{2}(x)h(y)+\lambda_{2}^{2}\,f(x)f(y),\;x,y\in S\end{equation}
for four unknown complex valued functions $f,g_{1},g_{2},h$ on a semigroup $S$ generated by its squares such that $f$ and $h$ are linearly independent, and $\lambda_{1},\lambda_{2}\in\mathbb{C}$ are given constants.\\
The system (\ref{eq01})-(\ref{eq02}) will play an important role in the resolution of a more general Levi-Civita functional equation: $f(xy)=g(x)h(y)+g(y)h(x)+k(x)k(y),\; x,y\in S.$
The solutions in \cite{C-K-N} and in \cite{Ajb and Elq1} are described in terms of multiplicative and additive functions. So, to solve a system of functional equations similar to (\ref{Eq8-0-0}), it is first necessary to give properties combining additive and multiplicative functions.\\
The outline of the paper is as follows.\\
The next section introduces some notations and terminology. In section 3, we prove some pertinent propositions that will be used in the prof of our main results. In section 4 we combine the authors's results \cite{Ajb and Elq1} about (\ref{Eq8-0-0}), with $f,h$ are linearly independent to get a more complete picture of the solutions of (\ref{Eq8-0-0}) in a semigroup generated by its squares. In section 5 we give the solutions of the system (\ref{eq01})-(\ref{eq02}) when $\lambda_{1}=0$ and $\lambda_{2}\neq0$. In section 6 the complete solutions of the system is given when $\lambda_{1}\neq0$ and $\lambda_{2}\neq0$.
\section{Notations and Terminology}
For a non empty set $X$ we denote by $\mathcal{F}(X,\mathbb{C})$ the $\mathbb{C}$-vector space of all maps from $X$ to $\mathbb{C}$.
Throughout this paper $S$ denotes a semigroup (a set with an
associative composition). We say that $S$ is generated by its squares if, for each $x\in S$ there exist $x_{1},\ldots{},x_{n}\in S$ such that
$x=x_{1}^{2}\cdot\cdot\cdot x_{n}^{2}$.\\
We say that a function $f:S\to\mathbb{C}$ is non-zero, and we denote $f\neq0$, if there exist an element $x\in S$ such that $f(x)\neq0$.\\
A additive function on $S$ is a homomorphic $a:S\to(\mathbb{C},+)$. We denote by $\mathcal{A}(S)$ the set of all additive functions on $S$.\\
A multiplicative function on $S$ is a homomorphic
$\chi:S\to(\mathbb{C},\cdot)$. We denote by $\mathcal{M}(S)$ the set of all multiplicative functions on $S$.\\If $\chi\in\mathcal{M}(S)$ and $\chi\neq0$, then
$I_{\chi}:=\{x\in S\mid\chi(x)=0\}$ is either empty or a proper
subset of $S$. $I_{\chi}$ is a two sided ideal in $S$ if not empty
and $S\setminus I_{\chi}$ is a subsemigroup of $S$.\\
For a function $g:S\to\mathbb{C}$ we define the set $$\mathcal{S}_{g}:=\{f:S\to\mathbb{C}\,|\,f(xy)=f(x)g(y)+g(x)f(y),\,x,y\in S\}.$$
For $\chi\in\mathcal{M}(S)$ and $\varphi\in\mathcal{F}(S\setminus I_{\chi},\mathbb{C})$ we define the function $\Psi_{\chi\,\varphi}:S\to\mathbb{C}$ by:
\begin{center}
\[\Psi_{\chi\,\varphi}(x)=
\left\{
\begin{array}{r c l}
\chi(x)\,\varphi(x)\quad\text{if}\quad x\in S\setminus I_{\chi}&\\0\quad\quad\quad\text{if}\quad x\in I_{\chi}&
\end{array}
\right.
\]
\end{center}
If $\chi\in\mathcal{M}(S)$ and $\chi\neq0$ we define the mapping $\Psi_{\chi}:\mathcal{F}(S\setminus I_{\chi},\mathbb{C})\to\mathcal{F}(S,\mathbb{C})$ by $\Psi_{\chi}(\varphi)=\Psi_{\chi\,\varphi}$ for all $\varphi\in\mathcal{F}(S\setminus I_{\chi},\mathbb{C})$.
\section{Auxiliary results}
In this section, we give some useful properties.
\begin{lem}\label{lem8-0}Let $T$ be a subemigroup of $S$. Let $a\in\mathcal{A}(T)$. Let $I$ be an ideal of $S$ such that $T\cap I\neq\emptyset$.
\par If $a(x)=0$ for all $x\in T\cap I$ then $a=0$.
\end{lem}
\begin{proof} Let $(x,y)\in T\times(T\cap I)$ be arbitrary. Since $T$ is a subemigroup of $S$ and $I$ is an ideal of $S$ we get that $xy\in T\cap I$. Ten $a(xy)=0$, which implies that $a(x)+a(y)=0$. As $a=0$ on $T\cap I$ by assumption, we obtain $a(x)=0$. So, $x\in T$ being arbitrary, we deduce that $a=0$.
\end{proof}
In Lemma \ref{lem8-2} we generalize \cite[Lemma 4.4]{Ajb and Elq1}.
\begin{lem}\label{lem8-2}Let $\mu_{1},\ldots{},\mu_{N}\in\mathcal{M}(S)$ be $N$ different multiplicative functions on $S$ and $A_{1},\ldots{},A_{N}\in\mathcal{A}(S)$.
\par If $\Sigma_{i=1}^{N}\mu_{i}\,A_{i}=\Sigma_{j=1}^{M}c_{j}\chi_{j}$, where
$c_{j}\in\mathbb{C}$ and $\chi_{j}\in\mathcal{M}(S)$ for each $j=1,2,\ldots{},N$, then $\mu_{i}\,A_{i}=0$ for each $i=1,2,\ldots{},N$.
\end{lem}
\begin{proof} We prove Lemma \ref{lem8-2} by induction on $N$. For $N=1$ we get Lemma \cite[Lemma 4.4]{Ajb and Elq1}. To get induction set up we assume that Lemma \ref{lem8-2} is true for an $N\geq1$ and prove it for $N+1$. So,
\begin{equation}\label{Eq8-0}\Sigma_{i=1}^{N+1}\mu_{i}\,A_{i}=\Sigma_{j=1}^{M}c_{j}\chi_{j}.\end{equation}
Let $x,y\in S$ be arbitrary. From the identity above we get that
\begin{equation}\label{Eq8-2}\Sigma_{i=1}^{N+1}\mu_{i}(xy)\,A_{i}(xy)=\Sigma_{j=1}^{M}c_{j}\chi_{j}(xy).\end{equation}
Using that $A_{1},\ldots{},A_{N}$ are additive and  $\mu_{1},\ldots{},\mu_{N},\chi_{1},\ldots{},\chi_{M}$ are multiplicative we get from (\ref{Eq8-2}) that
$$\Sigma_{i=1}^{N+1}\mu_{i}(x)\mu_{i}(y)\,A_{i}(x)+\Sigma_{i=1}^{N+1}\mu_{i}(x)\mu_{i}(y)\,A_{i}(y)=\Sigma_{j=1}^{M}c_{j}\chi_{j}(x)\chi_{j}(y).$$
Hence,
\begin{equation}\label{Eq8-3}\begin{split}&\Sigma_{i=1}^{N}\mu_{i}(x)\mu_{i}(y)\,A_{i}(x)+\mu_{N+1}(y)\mu_{N+1}(x)A_{N+1}(x)\\
&=\Sigma_{j=1}^{M}c_{j}\chi_{j}(x)\chi_{j}(y)-\Sigma_{i=1}^{N+1}\mu_{i}(y)\,A_{i}(y)\mu_{i}(x).\end{split}\end{equation}
As $\Sigma_{i=1}^{N+1}\mu_{i}\,A_{i}=\Sigma_{j=1}^{M}c_{j}\chi_{j}$ we get that \begin{equation}\label{Eq8-4}\mu_{N+1}(x)A_{N+1}(x)=\Sigma_{j=1}^{M}c_{j}\chi_{j}(x)-\Sigma_{i=1}^{N}\mu_{i}(x)\,A_{i}(x).\end{equation}
So, $x$ being arbitrary we get, by a small computation, from (\ref{Eq8-3}) and (\ref{Eq8-4}) that
\begin{equation}\label{Eq8-5}\begin{split}&\Sigma_{i=1}^{N}\mu_{i}((\mu_{i}(y)-\mu_{N+1}(y))\,A_{i})\\
&=\Sigma_{j=1}^{M}c_{j}(\chi_{j}(y)-\mu_{N+1}(y))\chi_{j}-\Sigma_{i=1}^{N+1}\mu_{i}(y)\,A_{i}(y)\mu_{i}.\end{split}\end{equation}
For each fixed $y\in S$ and $i=1,\ldots{},N$ the map $(\mu_{i}(y)-\mu_{N}(y))\,A_{i}:S\to\mathbb{C}$ is additive and $c_{j}(\chi_{j}(y)-\mu_{N+1}(y)),\,\mu_{i}(y)\,A_{i}(y)\in\mathbb{C}$, so by our induction hypothesis we deduce from (\ref{Eq8-5}) that $(\mu_{i}(y)-\mu_{N}(y))\mu_{i}\,A_{i}=0$ for each $y\in S$ and $i=1,\ldots{},N$. Since $\mu_{i}\neq\mu_{N}$ for each $i=1,\ldots{},N-1$, we deduce that $\mu_{i}\,A_{i}=0$ for all $i=1,\ldots{},N-1$. So, (\ref{Eq8-0}) reduces to $\mu_{N}\,A_{N}+\mu_{N+1}\,A_{N+1}=\Sigma_{j=1}^{M}c_{j}\chi_{j}$, which implies, by applying the induction assumption  that $\mu_{N+1}\,A_{N+1}=\mu_{N}\,A_{N}=0$. This completes the proof of Lemma \ref{lem8-2}.
\end{proof}
\begin{rem}\label{rem1} Let $\mu_{1},\mu_{2}\in\mathcal{M}(S)$ and $A_{1},A_{2}\in\mathcal{A}(S)$. If $\mu_{1}\,A_{1}+\mu_{2}\,A_{2}=\Sigma_{j=1}^{M}c_{j}\chi_{j}$, where
$c_{j}\in\mathbb{C}$ and $\chi_{j}\in\mathcal{M}(S)$ for each $j=1,2,\ldots{},N$, then $\mu_{1}\,A_{1}+\mu_{2}\,A_{2}=0$. Indeed,
\par If $\mu_{1}=\mu_{2}$ we get that $\mu_{1}\,(A_{1}+A_{2})=0$ by applying Lemma \cite[Lemma 4.4]{Ajb and Elq1}.
\par If $\mu_{1}\neq\mu_{2}$ we get that $\mu_{1}\,A_{1}=\mu_{1}\,A_{1}=0$ by applying Lemma \ref{lem8-2}.
\end{rem}
\begin{prop}\label{prop8-1}Let $\mu,\mu_{1},\mu_{2}\in\mathcal{M}(S)$ be non-zero multiplicative functions
on $S$ such that $\mu_{1}\neq\mu_{2}$, and $A\in\mathcal{A}(S\setminus I_{\mu}),A_{1}\in\mathcal{A}(S\setminus I_{\mu_{1}}),A_{2}\in\mathcal{A}(S\setminus I_{\mu_{2}})$.
\begin{enumerate}
\item\label{(1)} $\Psi_{\mu}$ is a injective linear map from $\mathcal{F}(S\setminus I_{\mu},\mathbb{C})$ to $ \mathcal{F}(S,\mathbb{C})$.
\item\label{(2)} Let $T$ be a subsemigroup of $S$. If $\Psi_{\mu}(A)=\Sigma_{j=1}^{N}c_{j}\chi_{j}$ on $T$, where
$c_{j}\in\mathbb{C}$ and $\chi_{j}\in\mathcal{M}(T)$ for each $j=1,2,\ldots{},N$, then $\Psi_{\mu}(A)=0$ on $T$.
\item\label{(3)} If $\Psi_{\mu_{1}}(A_{1})+\Psi_{\mu_{2}}(A_{2})=\Sigma_{j=1}^{N}c_{j}\chi_{j}$, where
$c_{j}\in\mathbb{C}$ and $\chi_{j}\in\mathcal{M}(S)$ for each $j=1,2,\ldots{},N$, then $\Psi_{\mu_{1}}(A_{1})+\Psi_{\mu_{2}}(A_{2})=0$.
\item\label{(4)} Let $T$ be a subsemigroup of $S$ such that $(S\setminus I_{\mu})\cap T\neq\emptyset$. Let $a,a_{1}\in\mathcal{A}(S\setminus I_{\mu})$. If $\Psi_{\mu}(a_{1}+a^{2})=\Sigma_{j=1}^{N}c_{j}\chi_{j}$ on $(S\setminus I_{\mu})\cap T$, where
$c_{j}\in\mathbb{C}$ and $\chi_{j}\in\mathcal{M}(S)$ for each $j=1,2,\ldots{},N$, then $a_{1}=0$ and $a=0$ on $(S\setminus I_{\mu})\cap T$.
\item\label{(5)} Let $a,a_{1}\in\mathcal{A}(S\setminus I_{\mu_{1}})$. If $\Psi_{\mu_{1}}(a_{1}+a^{2})+\Psi_{\mu}(A)=\Sigma_{j=1}^{N}c_{j}\chi_{j}$, where $c_{j}\in\mathbb{C}$ and $\chi_{j}\in\mathcal{M}(S)$ for each $j=1,2,\ldots{},N$, then $a=0$.
\item\label{(6)} Let $a,a_{1}\in\mathcal{A}(S\setminus I_{\mu_{1}})$. If $\Psi_{\mu}(A)=\Psi_{\mu_{1}}(a_{1}+a^{2})$, then $a=0$.
\item\label{(7)} Let $a,a_{1}\in\mathcal{A}(S\setminus I_{\mu_{1}})$. Assume that $a\neq0$ and $A\neq0$. If $\Psi_{\mu_{1}}(a_{1}+a^{2})=\Psi_{\mu}(A_{1}+A^{2})$, then $\mu=\mu_{1}$, $a_{1}=A_{1}$ and ($a=A$ or $a=-A$).
\end{enumerate}
\end{prop}
\begin{proof}(1) Let $\lambda\in\mathbb{C}$ and $\varphi_{1},\varphi_{2}\in\mathcal{F}(S\setminus I_{\mu},\mathbb{C})$. For all $x\in S\setminus I_{\mu}$ we have $\Psi_{\mu}(\lambda\varphi_{1})(x)=\mu(x)\,(\lambda\varphi_{1})(x)=\lambda(\mu(x)\,\varphi_{1}(x))=\lambda\Psi_{\mu}(\varphi_{1})(x)$ and $\Psi_{\mu}(\varphi_{1}+\varphi_{2})(x)=\mu(x)\,(\varphi_{1}+\varphi_{2})(x)=\mu(x)\,\varphi_{1}(x)+\mu(x)\,\varphi_{2}(x)
=\Psi_{\mu}(\varphi_{1})(x)+\Psi_{\mu}(\varphi_{2})(x)=(\Psi_{\mu}(\varphi_{1})+\Psi_{\mu}(\varphi_{2}))(x)$.
By a similar computation, we check that $\Psi_{\mu}(\lambda\varphi_{1})(x)=\lambda\Psi_{\mu}(\varphi_{1})(x)$
 and $\Psi_{\mu}(\varphi_{1}+\varphi_{2})(x)=(\Psi_{\mu}(\varphi_{1})+\Psi_{\mu}(\varphi_{2}))(x)$ for all $x\in I_{\mu}$. Hence, $\Psi_{\mu}(\lambda\varphi_{1})=\lambda\Psi_{\mu}(\varphi_{1})$ and $\Psi_{\mu}(\varphi_{1}+\varphi_{2})=\Psi_{\mu}(\varphi_{1})+\Psi_{\mu}(\varphi_{2})$.\\
Moreover, if $\Psi_{\mu}(\varphi_{1})=0$ then $\mu(x)\,\varphi_{1}(x)=0$ for all $x\in S\setminus I_{\mu}$, which implies that $\varphi_{1}=0$ because $\mu(x)\neq0$  for all $x\in S\setminus I_{\mu}$. Hence, $\Psi_{\mu}$ is a injective linear map from $\mathcal{F}(S\setminus I_{\mu},\mathbb{C})$ to $ \mathcal{F}(S,\mathbb{C})$.\\
(2) If $(S\setminus I_{\mu})\cap T=\emptyset$ or $I_{\mu}\cap T=\emptyset$ then $T\subset I_{\mu}$ or $T\subset S\setminus I_{\mu}$. Hence, $\Psi_{\mu}(A)=0$ on $T$ or $\mu\,A=\Sigma_{j=1}^{N}c_{j}\chi_{j}$ on $T$. We deduce, by applying Lemma \cite[Lemma 4.4]{Ajb and Elq1} on the subsemigroup $T$ to the last identity, that $\Psi_{\mu}(A)=0$ on $T$.\\
Now, we assume that $(S\setminus I_{\mu})\cap T\neq\emptyset$ and $I_{\mu}\cap T\neq\emptyset$. Then $\mu\,A=\Sigma_{j=1}^{N}c_{j}\chi_{j}$ on $(S\setminus I_{\mu})\cap T$. By applying Lemma \cite[Lemma 4.4]{Ajb and Elq1} on the subsemigroup $(S\setminus I_{\mu})\cap T$ we get that $\Psi_{\mu}(A)=0$ on $(S\setminus I_{\mu})\cap T$. Moreover, we have $\Psi_{\mu}(A)=0$ on $I_{\mu}\cap T$. We infer that $\Psi_{\mu}(A)=0$.\\
(3) If $I_{\mu_{1}}=I_{\mu_{2}}=\emptyset$ then $\Psi_{\mu_{1}}(A_{1})=\mu_{1}\,A_{1}$ and $\Psi_{\mu_{2}}(A_{2})=\mu_{2}\,A_{2}$. Hence, $\mu_{1}\,A_{1}+\mu_{2}\,A_{2}=\Sigma_{j=1}^{N}c_{j}\chi_{j}$. So, according to Remark \ref{rem1}, we deduce that $\Psi_{\mu_{1}}(A_{1})+\Psi_{\mu_{2}}(A_{2})=0$.\\ Now, assume that $I_{\mu_{1}}\neq\emptyset$. Then, $\Psi_{\mu_{1}}(A_{1})=0$ on $I_{\mu_{1}}$. So, $\Psi_{\mu_{2}}(A_{2})=\Sigma_{j=1}^{N}c_{j}\chi_{j}$ on the subsemigroup $I_{\mu_{1}}$. Hence,  by applying Proposition \ref{prop8-1}(\ref{(2)}) on the subsemigroup $I_{\mu_{1}}$, we obtain $\Psi_{\mu_{2}}(A_{2})=0$ on $I_{\mu_{1}}$. So that \begin{equation}\label{Eq8-6-1}\Psi_{\mu_{1}}(A_{1})+\Psi_{\mu_{2}}(A_{2})=0\,\,\text{on}\,\,I_{\mu_{1}}.\end{equation} We consider the following cases:\\
\underline{Case 1}: $(S\setminus I_{\mu_{1}})\cap (S\setminus I_{\mu_{2}})\neq\emptyset$. Then
\begin{equation}\label{Eq8-6}\Psi_{\mu_{1}}(A_{1})+\Psi_{\mu_{2}}(A_{2})=0\,\,\text{on}\,\,(S\setminus I_{\mu_{1}})\cap (S\setminus I_{\mu_{2}}).\end{equation}
Indeed,
\begin{equation}\label{Eq8-7}\mu_{1}\,A_{1}+\mu_{2}\,A_{2}=\Sigma_{j=1}^{N}c_{j}\chi_{j}\end{equation}
on the subsemigroup $(S\setminus I_{\mu_{1}})\cap (S\setminus I_{\mu_{2}})$. So, according to Remark \ref{rem1} and seeing that $\Psi_{\mu_{1}}(A_{1})=\mu_{1}\,A_{1}$ and $\Psi_{\mu_{2}}(A_{2})=\mu_{2}\,A_{2}$ on $(S\setminus I_{\mu_{1}})\cap (S\setminus I_{\mu_{2}})$,
we get (\ref{Eq8-6}).\\
\underline{Subcase 1.1}: $(S\setminus I_{\mu_{1}})\cap I_{\mu_{2}}\neq\emptyset$. Since $\Psi_{\mu_{2}}(A_{2})=0$ on $I_{\mu_{2}}$ we get that $\Psi_{\mu_{1}}(A_{1})=\Sigma_{j=1}^{N}c_{j}\chi_{j}$ on the subsemigroup $(S\setminus I_{\mu_{1}})\cap I_{\mu_{2}}$. So, by applying Proposition \ref{prop8-1}(\ref{(2)}) on the subsemigroup $(S\setminus I_{\mu_{1}})\cap I_{\mu_{2}}$, we get that $\Psi_{\mu_{1}}(A_{1})=0$ on $(S\setminus I_{\mu_{1}})\cap I_{\mu_{2}}$. It follows that $\Psi_{\mu_{1}}(A_{1})+\Psi_{\mu_{2}}(A_{2})=0$ on $(S\setminus I_{\mu_{1}})\cap I_{\mu_{2}}$. Combining this with (\ref{Eq8-6}) we deduce that
\begin{equation*}\Psi_{\mu_{1}}(A_{1})+\Psi_{\mu_{2}}(A_{2})=0\,\,\text{on}\,\,S\setminus I_{\mu_{1}}.\end{equation*}
 So, taking (\ref{Eq8-6-1}) into account, we get $\Psi_{\mu_{1}}(A_{1})+\Psi_{\mu_{2}}(A_{2})=0$.\\
\underline{Subcase 1.2}: $(S\setminus I_{\mu_{1}})\cap I_{\mu_{2}}=\emptyset$. Then, $S\setminus I_{\mu_{1}}\subset S\setminus I_{\mu_{2}}$. So, (\ref{Eq8-6}) implies that $\Psi_{\mu_{1}}(A_{1})+\Psi_{\mu_{2}}(A_{2})=0\,\,\text{on}\,\,S\setminus I_{\mu_{1}}$. Hence, as in Subcase 1.1, we deduce that $\Psi_{\mu_{1}}(A_{1})+\Psi_{\mu_{2}}(A_{2})=0$.\\
\underline{Case 2}: $(S\setminus I_{\mu_{1}})\cap (S\setminus I_{\mu_{2}})=\emptyset$. Then $S\setminus I_{\mu_{1}}\subset I_{\mu_{2}}$ and $S\setminus I_{\mu_{2}}\subset I_{\mu_{1}}$.
As $\Psi_{\mu_{1}}(A_{1})+\Psi_{\mu_{2}}(A_{2})=\Sigma_{j=1}^{N}c_{j}\chi_{j}$,
$\Psi_{\mu_{1}}(A_{1})=0$ on $I_{\mu_{1}}$ and $\Psi_{\mu_{2}}(A_{2})=0$ on $I_{\mu_{2}}$ we get that $$\Psi_{\mu_{1}}(A_{1})=\Sigma_{j=1}^{N}c_{j}\chi_{j}\,\,\text{on}\,\,S\setminus I_{\mu_{1}}$$ and
$$\Psi_{\mu_{2}}(A_{2})=\Sigma_{j=1}^{N}c_{j}\chi_{j}\,\,\text{on}\,\,S\setminus I_{\mu_{2}}.$$
It follows, according to Proposition \ref{prop8-1}(\ref{(2)}), that $\Psi_{\mu_{1}}(A_{1})=0$ and $\Psi_{\mu_{2}}(A_{2})=0$. Then $\Psi_{\mu_{1}}(A_{1})+\Psi_{\mu_{2}}(A_{2})=0$.\\
(4) On the subsemigroup $(S\setminus I_{\mu})\cap T$ we get that $\mu\,a_{1}+\mu\,a^{2}=\Sigma_{j=1}^{N}c_{j}\chi_{j}$.  We can divide by $\mu(x)$ for all $x\in (S\setminus I_{\mu})\cap T$. So, by putting $m_{j}(x):=\dfrac{\chi_{j}(x)}{\mu(x)}$ for all $x\in (S\setminus I_{\mu})\cap T$ and for each $j=1,2,\ldots{},N$, we get that
\begin{equation}\label{Eq8-8}a_{1}+a^{2}=\Sigma_{j=1}^{N}c_{j}m_{j}\,\,\text{on}\,\,(S\setminus I_{\mu})\cap T.\end{equation} Let $x,y\in (S\setminus I_{\mu})\cap T$ be arbitrary. Since $a,a_{1}:S\setminus I_{\mu}\to\mathbb{C}$ are  additive and $m_{j}:S\setminus I_{\mu}\to\mathbb{C}$ is multiplicative for each $j=1,2,\ldots{},N$, we get from (\ref{Eq8-8}) that $$a_{1}(x)+a_{1}(y)+a^{2}(x)+a^{2}(y)+2a(x)a(y)=\Sigma_{j=1}^{N}c_{j}m_{j}(x)m_{j}(y),$$ which implies that $$\Sigma_{j=1}^{N}c_{j}m_{j}(x)+\Sigma_{j=1}^{N}c_{j}m_{j}(y)+2a(x)a(y)=\Sigma_{j=1}^{N}c_{j}m_{j}(x)m_{j}(y).$$ Hence, $x$ and $y$ being arbitrary, we deduce from the identity above that
\begin{equation}\label{Eq8-10}1\cdot(2a(y)a)=\Sigma_{j=1}^{N}c_{j}(m_{j}(y)-1)m_{j}-(\Sigma_{j=1}^{N}c_{j}m_{j}(y))\cdot1\,\,\text{on}\,\,(S\setminus I_{\mu})\cap T.\end{equation}
So, applying Lemma \cite[Lemma 4.4]{Ajb and Elq1} we get that $a(y)a=0$ for all $y\in (S\setminus I_{\mu})\cap T$, then $a=0$ on $(S\setminus I_{\mu})\cap T$. It follows that (\ref{Eq8-8}) reduces to $a_{1}=\Sigma_{j=1}^{N}c_{j}m_{j}$ on $(S\setminus I_{\mu})\cap T$. So, according to Lemma \cite[Lemma 4.4]{Ajb and Elq1}, we get that $a_{1}=0$ on $(S\setminus I_{\mu})\cap T$.\\
(5) If $(S\setminus I_{\mu_{1}})\cap(S\setminus I_{\mu})=\emptyset$ then $S\setminus I_{\mu_{1}}\subset I_{\mu}$. Hence, $\Psi_{\mu}(A)=0$ on $S\setminus I_{\mu_{1}}$. So that, $\Psi_{\mu_{1}}(a_{1}+a^{2})=\Sigma_{j=1}^{N}c_{j}\chi_{j}$ on the subsemigroup $S\setminus I_{\mu_{1}}$. Thus, by applying Proposition \ref{prop8-1}(\ref{(4)}), we derive that $a=0$. So, we assume that $(S\setminus I_{\mu_{1}})\cap(S\setminus I_{\mu})\neq\emptyset$. Then, $\mu_{1}(a_{1}+a^{2})+\mu\,A=\Sigma_{j=1}^{N}c_{j}\chi_{j}$. Let $m:=\dfrac{\mu}{\mu_{1}}$ and $m_{j}:=\dfrac{\chi_{j}}{\mu_{1}}$ on $(S\setminus I_{\mu_{1}})\cap(S\setminus I_{\mu})$ for each $j=1,2,\ldots{},N$. When we divide the last identity by $\mu_{1}$ we obtain
\begin{equation}\label{Eq8-10-1}a_{1}+a^{2}+m\,A=\Sigma_{j=1}^{N}c_{j}m_{j}.\end{equation}
Since $a_{1},a,A$ are additive and $m,m_{j}$ are multiplicative for each $j=1,2,\ldots{},N$, we derive from (\ref{Eq8-10-1}) that \begin{equation*}\begin{split}&a_{1}(x)+a^{2}(x)+a_{1}(y)+a^{2}(y)+2a(y)a(x)+m(y)m(x)A(x)+m(y)A(y)m(x)\\
&\quad\quad\quad\quad\quad=\Sigma_{j=1}^{N}c_{j}m_{j}(y)m_{j}(x)\end{split}\end{equation*}
for all $x,y\in(S\setminus I_{\mu_{1}})\cap(S\setminus I_{\mu})$. Hence, taking (\ref{Eq8-10-1}) into account, we get that
\begin{equation*}\begin{split}&\Sigma_{j=1}^{N}c_{j}m_{j}-m\,A+a_{1}(y)+a^{2}(y)+2a(y)a+m(y)m\,A+m(y)A(y)m\\
&\quad\quad\quad\quad\quad=\Sigma_{j=1}^{N}c_{j}m_{j}(y)m_{j}\end{split}\end{equation*} for all $y\in(S\setminus I_{\mu_{1}})\cap(S\setminus I_{\mu})$. So that
\begin{equation}\label{Eq8-10-2}m\,((m(y)-1)A)+1\cdot(2a(y)a)=\Sigma_{j=1}^{N}c_{j}(m_{j}(y)-1)m_{j}-m(y)A(y)m-(a_{1}(y)+a^{2}(y))\cdot1\end{equation} for all $y\in(S\setminus I_{\mu_{1}})\cap(S\setminus I_{\mu})$. Then
\begin{equation}\label{Eq8-10-3}a=0\,\,\text{on}\,\,(S\setminus I_{\mu_{1}})\cap(S\setminus I_{\mu}).\end{equation}
Indeed,
\par if $m=1$ then, in view of (\ref{Eq8-10-1}), the identity (\ref{Eq8-10-2}) reduces to $$1\cdot(2a(y)a)=\Sigma_{j=1}^{N}c_{j}(m_{j}(y)-1)m_{j}-(\Sigma_{j=1}^{N}c_{j}m_{j}(y))\cdot1,$$
for all $y\in(S\setminus I_{\mu_{1}})\cap(S\setminus I_{\mu})$, which implies, according to Lemma \cite[Lemma 4.4]{Ajb and Elq1}, that $a(y)a=0$ on $(S\setminus I_{\mu_{1}})\cap(S\setminus I_{\mu})$ for all $y\in(S\setminus I_{\mu_{1}})\cap(S\setminus I_{\mu})$. Hence, we get (\ref{Eq8-10-3}).
\par If $m\neq1$ then, by applying Lemma \ref{lem8-2}, we get from (\ref{Eq8-10-2}) that $a(y)a=0$  on $(S\setminus I_{\mu_{1}})\cap(S\setminus I_{\mu})$ for all $y\in(S\setminus I_{\mu_{1}})\cap(S\setminus I_{\mu})$. So, we obtain (\ref{Eq8-10-3}).\\
Now, we consider the following cases:\\
\underline{Case 1}: $(S\setminus I_{\mu_{1}})\cap I_{\mu}=\emptyset$. Then $(S\setminus I_{\mu_{1}})\subset(S\setminus I_{\mu})$, which implies, taking (\ref{Eq8-10-3}) into account, that $a=0$.\\
\underline{Case 2}: $(S\setminus I_{\mu_{1}})\cap I_{\mu}\neq\emptyset$. Since $\Psi_{\mu}(A)=0$ on $I_{\mu}$ we get that $\Psi_{\mu_{1}}(a_{1}+a^{2})=\Sigma_{j=1}^{N}c_{j}\chi_{j}$ on the subsemigroup $(S\setminus I_{\mu_{1}})\cap I_{\mu}$. Hence, by applying Proposition \ref{prop8-1}(\ref{(4)}) for $T=I_{\mu}$, we get that $a=0$ on $(S\setminus I_{\mu_{1}})\cap I_{\mu}$. So, taking (\ref{Eq8-10-3}) into account, we obtain $a=0$.\\
(6) We derive the result from (5) by writing $-A$ instead of $A$ and taking $c_{j}=0$ for $j=1,2,\ldots{},N$.\\
(7) If $(S\setminus I_{\mu})\cap I_{\mu_{1}}\neq\emptyset$ then we get from the identity $\Psi_{\mu_{1}}(a_{1}+a^{2})=\Psi_{\mu}(A_{1}+A^{2})$ that $\Psi_{\mu}(A_{1}+A^{2})=0$ on $(S\setminus I_{\mu})\cap I_{\mu_{1}}$. So, by applying Proposition \ref{prop8-1}(\ref{(1)}) on the semigroup $(S\setminus I_{\mu})\cap I_{\mu_{1}}$, we deduce that $A_{1}+A^{2}=0$ on $(S\setminus I_{\mu})\cap I_{\mu_{1}}$. So, like in Subcase 1.1 of part (6), we get that $A=0$ on $(S\setminus I_{\mu})\cap I_{\mu_{1}}$. Then, according to Lemma \ref{lem8-0}, we infer that $A=0$, which contradicts the assumption on $A$. Hence, $(S\setminus I_{\mu})\cap I_{\mu_{1}}=\emptyset$, which implies that $S\setminus I_{\mu}\subset S\setminus I_{\mu_{1}}$. Similarly, we prove that $S\setminus I_{\mu_{1}}\subset S\setminus I_{\mu}$. So that $S\setminus I_{\mu}=S\setminus I_{\mu_{1}}$. Let $\dfrac{\mu}{\mu_{1}}=:m\in\mathcal{M}(S\setminus I_{\mu})$. Then the identity $\Psi_{\mu_{1}}(a_{1}+a^{2})=\Psi_{\mu}(A_{1}+A^{2})$ reduces to
\begin{equation*}a_{1}+a^{2}=m(A_{1}+A^{2})\end{equation*}
on $S\setminus I_{\mu}$. Using similar computations to the ones in part (5) we obtain from the identity above that
\begin{equation}\label{Eq8-14}\begin{split}&m[(1-m(y))A_{1}-2m(y)A(y)A+(1-m(y))A^{2}]+1\cdot(2a(y)a)\\
&=-(a_{1}(y)+a^{2}(y))\cdot1+m(y)(A_{1}(y)+A^{2}(y))m\end{split}\end{equation}
for each $y\in S\setminus I_{\mu}$.
\par If $m\neq1$ then there exists $y_{0}\in S\setminus I_{\mu}$ such that $1-m(y_{0})\neq0$. Let $1-m(y_{0})=:\lambda^{2}$. As $(1-m(y_{0}))A_{1}-2m(y_{0})A(y_{0})A,\lambda\,A\in\mathcal{A}(S\setminus I_{\mu})$, we derive from (\ref{Eq8-14}), by applying Proposition \ref{prop8-1}(\ref{(5)}) on the subsemigroup $S\setminus I_{\mu}$, that $\lambda\,A=0$. Seeing that $\lambda\neq0$, we get that $A=0$. Contradicting the assumption on $A$.
\par Hence $m=1$. So, $\mu=\mu_{1}$. It follows that $\Psi_{\mu}(a_{1}+a^{2})=\Psi_{\mu}(A_{1}+A^{2})$. Then, according to Proposition \ref{prop8-1}(\ref{(1)}), we get that
\begin{equation}\label{Eq8-15}a_{1}+a^{2}=A_{1}+A^{2}.\end{equation}
So, $a^{2}-A^{2}\in\mathcal{A}(S\setminus I_{\mu})$, which implies that $$a^{2}(xy)-A^{2}(xy)=a^{2}(x)-A^{2}(x)+a^{2}(y)-A^{2}(y)$$ for all $x,y\in S\setminus I_{\mu}$. Since $a,A\in \mathcal{A}(S\setminus I_{\mu})$,
we derive from the identity above, by a small computation, that
\begin{equation}\label{Eq8-16}a(x)a(y)=A(x)A(y)\end{equation} for all $x,y\in S\setminus I_{\mu}$. Then, $a^{2}=A^{2}$ which implies, taking (\ref{Eq8-15}) into account, that $a_{1}=A_{1}$. Moreover, seeing that $a\neq0$, we derive from (\ref{Eq8-16}) that there exists a constant $\alpha\in\mathbb{C}\setminus\{0\}$ such that $a=\alpha\,A$. Since $a^{2}=A^{2}$, we get that $\alpha^{2}=1$. Hence, $\alpha=1$ or $\alpha=-1$. So that, $a=A$ or $a=-A$.
\end{proof}
\begin{prop}\label{prop8-1-1} Let $S$ be a semigroup and $g,h:S\to\mathbb{C}$ be functions. Then
\begin{enumerate}
\item\label{(1-1)} $\mathcal{S}_{g}$ is a $\mathbb{C}-$ vector space.
\item\label{(1-2)} If there exists $f\neq0\in\mathcal{S}_{g}\cap\mathcal{S}_{h}$ then $g=h$.
\item Let $\mu,\chi\in\mathcal{M}(S)\setminus\{0\},\,a\in\mathcal{A}(S\setminus I_{\chi})\setminus\{0\}$ or $A\in\mathcal{A}(S\setminus I_{\mu})\setminus\{0\}$. If $\Psi_{\chi}(a)=\Psi_{\mu}(A)$ then $\mu=\chi$ and $a=A$.
\end{enumerate}
\end{prop}
\begin{proof} (1) Notice that $\mathcal{S}_{g}$ is a subset of the $\mathbb{C}-$ vector space $\mathcal{F}(X,\mathbb{C})$. We have $\mathcal{S}_{g}\neq\emptyset$ because the null function belongs to $\mathcal{S}_{g}$. Moreover, if $f_{1},f_{2}\in\mathcal{S}_{g}$ and $\alpha\in\mathbb{C}$ then
\begin{equation*}\begin{split}&(\alpha\,f_{1}+f_{2})(xy)=\alpha\,f_{1}(xy)+f_{2}(xy)\\
&\quad\quad\quad\quad\quad\quad\,\,\,=\alpha\,f_{1}(x)g(y)+\alpha\,f_{1}(y)g(x)+f_{2}(x)g(y)+f_{2}(y)g(x)\\
&\quad\quad\quad\quad\quad\quad\,\,\,=(\alpha\,f_{1}+f_{2})(x)g(y)+g(x)(\alpha\,f_{1}+f_{2})(y),\end{split}\end{equation*} for all $x,y\in S$. So, $\alpha\,f_{1}+f_{2}\in \mathcal{S}_{g}$. Hence, $\mathcal{S}_{g}$ is a vector subspace of the $\mathbb{C}-$ vector space $\mathcal{F}(X,\mathbb{C})$. We deduce that $\mathcal{S}_{g}$ is a $\mathbb{C}-$ vector space.\\
(2) We have $f(xy)=f(x)g(y)+g(x)f(y)=f(x)h(y)+h(x)f(y)$ for all $x,y\in S$. Then
\begin{equation}\label{eq}f(x)(g(y)-h(y))=-f(y)(g(x)-h(x))\end{equation} for all $x,y\in S$. As $f\neq 0$  there exists $y_{0}\in S$ such that $f(y_{0})\neq0$. So, (\ref{eq}) implies that there exists $\alpha\in \mathbb{C}$ such that $g-h=\alpha\,f$. Substituting this in (\ref{eq}) we get that $\alpha\,f(x)f(y)=0$ for all $x,y\in S$. By putting $x=y=y_{0}$ in this identity we infer that $\alpha=0$. Hence, $g=h$.\\
(3) Assume that $a\in\mathcal{A}(S\setminus I_{\chi})\setminus\{0\}$. Let $x,y\in S$.
\par If $x,y\in S\setminus I_{\chi}$ then $xy\in S\setminus I_{\chi}$. So, $\Psi_{\chi}(a)(xy)=\chi(xy)a(xy)=\chi(x)a(x)\chi(y)+\chi(x)\chi(y)a(y)=\Psi_{\chi}(a)(x)\chi(y)+\chi(x)\Psi_{\chi}(a)(y)$.
\par If $x\in I_{\chi}$ or $y\in I_{\chi}$ then $xy\in I_{\chi}$. So that $\Psi_{\chi}(a)(xy)=0$ and, ($\chi(x)=0$ and $\Psi_{\chi}(a)(x)=0$) or ($\chi(y)=0$ and $\Psi_{\chi}(a)(y)=0$). Hence, $\Psi_{\chi}(a)(x)\chi(y)+\chi(x)\Psi_{\chi}(a)(y)=0=\Psi_{\chi}(a)(xy)$.
\par So, $\Psi_{\chi}(a)(xy)=\Psi_{\chi}(a)(x)\chi(y)+\chi(x)\Psi_{\chi}(a)(y)$ for all $x,y\in S$. Then $\Psi_{\chi}(a)\in\mathcal{S}_{\chi}$. Similarly we have $\Psi_{\mu}(A)\in\mathcal{S}_{\mu}$. Then $\Psi_{\chi}(a)\in\mathcal{S}_{\chi}\cap\mathcal{S}_{\mu}$. Moreover, since $a\neq0$ we get, according to Proposition \ref{prop8-1}(1), that $\Psi_{\chi}(a)\neq0$. So, applying Proposition \ref{prop8-1-1}(2), we deduce that $\chi=\mu$ and then, using Proposition \ref{prop8-1}(1), $a=A$.
\end{proof}
\section{Preliminaries}
Throughout the rest of this article we require our  semigroup $S$ to be generated by its squares. First, we give the solutions of the sine addition law. According to \cite[Proposition 3.1]{STETKAER} we have the following
\begin{prop}\label{prop8-1-2}The solutions $f,g:S\to\mathbb{C}$ of the functional equation
\begin{equation*}f(xy)=f(x)g(y)+g(x)f(y),\,x,y\in S\end{equation*} such that $f\neq0$ are the following, where $\chi,\chi_{1},\chi_{2}\in\mathcal{M}(S)$  such that $\chi_{1}\neq\chi_{2}$ and $\chi\neq0$, $A\in\mathcal{A}(S\setminus I_{\chi})$ such that $A\neq0$ and $\alpha\in\mathbb{C}\setminus\{0\}$ is a constant.
\begin{enumerate}
\item $f=\alpha(\chi_{1}-\chi_{2})$, $g=\dfrac{\chi_{1}+\chi_{2}}{2}$.
\item $f=\Psi_{\chi}(A)$, $g=\chi$.
\end{enumerate}
\end{prop}
Next, we give a small improvement to the description \cite[Theorem 4.3]{Ajb and Elq1} of the solutions of the cosine-sine functional equation (\ref{Eq8-0-0}).
\begin{prop}\label{prop8-2}The solutions $f,g,h:S\to\mathbb{C}$ of the functional equation
\begin{equation*}f(xy)=f(x)g(y)+g(x)f(y)+h(x)h(y),\,x,y\in S,\end{equation*} such that $\{f,h\}$ is linearly independent, are of the form
$$f=F,\,\,\,\,g=-\dfrac{1}{2}\delta^{2}\,F+G+\delta\,H,\,\,\,\,h=-\delta
\,F+H,$$
where $\delta\in\mathbb{C}$ is a constant and $F,G,H:S\to\mathbb{C}$ are of the following families, where $\mu,\chi,\chi_{1},\chi_{2},\chi_{3}\in\mathcal{M}(S)$ such that $\chi\neq0$, $\mu\neq\chi$ and $\chi_{1},\chi_{2},\chi_{3}$ are different, $A,A_{1}\in\mathcal{A}(S\setminus I_{\chi})$ such that $A\neq0$ and $c,d,\alpha,\beta\in\mathbb{C}\setminus\{0\}$ are constants.
\begin{enumerate}
\item\label{(8)} $$F=\dfrac{1}{2}\Psi_{\chi}(A_{1}+A^{2}),\,\,\,\,G=\chi,\,\,\,\,H=\Psi_{\chi}(A).$$
\item\label{(9)} $$F=c^{2}(\mu-\chi)-c\,\Psi_{\chi}(A),\,\,\,G=\chi,\,\,\,\,H=c(\mu-\chi).$$
\item\label{(10)} $$F=c(\mu-\chi)+cd\,\Psi_{\chi}(A),\,\,\,\,G=\dfrac{\mu+\chi}{2}-\dfrac{1}{2}d\,\Psi_{\chi}(A),\,\,\,\,H=\Psi_{\chi}(A)$$ with $1-cd^{2}=0$.
\item\label{(11)} $$F=c\,\beta\,\chi_{1}+c\,(2-\beta)\,\chi_{2}-2\,c\,\chi_{3},$$
    $$G=\dfrac{1}{4}\,\beta\,\chi_{1}+\dfrac{1}{4}\,(2-\beta)\,\chi_{2}+\dfrac{1}{2}\,\chi_{3},$$$$H=\dfrac{\chi_{1}-\chi_{2}}{2\alpha}$$
    with $2\,c\,\alpha^{2}\,\beta\,(2-\beta)=1$.
\end{enumerate}
\end{prop}
\section{Solution of the system (\ref{eq01})-(\ref{eq02}) when $\lambda_{1}=0$ and $\lambda_{2}\neq0$}
In this section we solve the following system, which we denote $(\ref{EQ1})-(\ref{EQ2})$:
\begin{align}
\label{EQ1}f(xy)=f(x)g_{1}(y)+g_{1}(x)f(y),\;x,y\in S,\\
\label{EQ2}h(xy)=h(x)g_{2}(y)+g_{2}(x)h(y)+\lambda_{2}^{2}\,f(x)f(y),\;x,y\in S
\end{align}
for four unknown complex valued functions $f,g_{1},g_{2},h$ on a semigroup $S$ generated by its squares such that $f$ and $h$ are linearly independent, and $\lambda_{2}\in\mathbb{C}\setminus\{0\}$ is a given constant.\\
The solutions of the system $(\ref{EQ1})-(\ref{EQ2})$ are given in Theorem \ref{thm1}.
\begin{thm}\label{thm1} The solutions $f,g_{1},h,g_{2}:S\to \mathbb{C}$ of the system $(\ref{EQ1})-(\ref{EQ2})$ such that  $f$ and $h$ are linearly independent are of the following forms where $m,\mu,\chi_{1},\chi_{2},\chi_{3}\in\mathcal{M}(S)$ such that $m\neq0$, $m\neq\mu$ and $\chi_{1},\chi_{2},\chi_{3}$ are different, $A,A_{1}\in\mathcal{A}(S\setminus I_{m})$ such that $A\neq0$, and $\alpha,\beta,\lambda,c,d\setminus\{0\}$ are constants.
 \begin{enumerate}
\item\label{(t1)} $$f=\Psi_{m}(A),\,\,g_{1}=m,$$ $$h=\dfrac{1}{2}\Psi_{m}(A_{1}+\lambda_{2}^{2}\,A^{2}),\,\,g_{2}=m.$$
\item\label{(t2)} $$f=\alpha(\mu-m),\,\,g_{1}=\dfrac{\mu+m}{2},$$
                  $$h=\alpha^{2}\,\lambda_{2}^{2}(\mu-m)-\alpha\,\lambda_{2}\,\Psi_{m}(A),\,\,g_{2}=m.$$
\item\label{(t3)} $$f=\dfrac{1}{\lambda_{2}}\,\Psi_{m}(A),\,\,g_{1}=m,$$
$$h=c^{2}(\mu-m)-c\,\Psi_{m}(A),\,\,g_{2}=\dfrac{\mu-m}{2}-\dfrac{1}{2c}\,\Psi_{m}(A).$$
\item\label{(t4)} $$f=-\dfrac{1}{d\,\lambda_{2}}(\mu-m),\,\,g_{1}=\dfrac{\mu+m}{2},$$
                  $$h=c(\mu-m)+\dfrac{1}{d}\,\Psi_{m}(A),\,\,g_{2}=m.$$
\item\label{(t5)} $$f=\dfrac{1}{\lambda_{2}}\Psi_{m}(A),\,\,g_{1}=m,$$
$$h=c(\mu-m)+\dfrac{1}{d}\,\Psi_{m}(A),\,\,g_{2}=\dfrac{\mu+m}{2}-\dfrac{1}{2}d\,\Psi_{m}(A)$$
with $1-cd^{2}=0$.
\item\label{(t6)} $$f=\dfrac{\chi_{1}-\chi_{2}}{2\,\lambda\,\lambda_{2}},\,\,g_{1}=\dfrac{\chi_{1}+\chi_{2}}{2},$$
$$h=c\,\beta\,\chi_{1}+c\,(2-\beta)\,\chi_{2}-2\,c\,\chi_{3},\,\,g_{2}=\dfrac{\beta\,\chi_{1}+(2-\beta)\,\chi_{2}+2\,\chi_{3}}{4},$$
with $2\,c\,\lambda^{2}\,\beta\,(2-\beta)=1$.
\item\label{(t7)} $$f=\dfrac{1}{\lambda\,\lambda_{2}(2-\beta)}(\chi_{1}-\chi_{3}),\,\,g_{1}=\dfrac{\chi_{1}+\chi_{3}}{2},$$
$$h=c\,\beta\,\chi_{1}+c\,(2-\beta)\,\chi_{2}-2\,c\,\chi_{3},\,\,g_{2}=\dfrac{-\beta\,\chi_{1}+(2-\beta)\,\chi_{2}+2\,\chi_{3}}{2(2-\beta)},$$ with $2\,c\,\lambda^{2}\,\beta\,(2-\beta)=1$
\item\label{(t8)} $$f=-\dfrac{1}{\lambda\,\lambda_{2}\,\beta}(\chi_{2}-\chi_{3}),\,\,g_{1}=\dfrac{\chi_{2}+\chi_{3}}{2},$$
$$h=c\,\beta\,\chi_{1}+c\,(2-\beta)\,\chi_{2}-2\,c\,\chi_{3},\,\,g_{2}=\dfrac{\beta\,\chi_{1}-(2-\beta)\,\chi_{2}+2\,\chi_{3}}{2\,\beta},$$
 with $2\,c\,\lambda^{2}\,\beta\,(2-\beta)=1$.
\end{enumerate}
\end{thm}
\begin{proof}
From (\ref{EQ2}) we deduce,  according to Proposition \ref{prop8-2}, that
\begin{equation}\label{EQ3} h=H \end{equation}
\begin{equation}\label{EQ4} g_{2}=-\dfrac{1}{2}\delta^{2}\,H+G+\lambda_{2}\delta\,F \end{equation}
\begin{equation}\label{EQ5} \lambda_{2}\;f=-\delta\,H+\lambda_{2}\,F,\end{equation}
where $\delta\in\mathbb{C}$ is a constant and $\lambda_{2}\,F,G,H:S\to \mathbb{C}$ are of the forms in Proposition \ref{prop8-2}(\ref{(8)}),(\ref{(9)}),(\ref{(10)}),(\ref{(11)}) with the same constraints. So, we split the discussion into the following cases:\\
\underline{Case 1}: $H=\dfrac{1}{2}\Psi_{m}(a_{1}+a^{2}),\,\,\,\,G=m,\,\,\,\,\lambda_{2}\,F=\Psi_{m}(a)$, where $0\neq m\in\mathcal{M}(S)$, $a\neq0,a_{1}\in\mathcal{A}(S\setminus I_{m})$. Applying Proposition \ref{prop8-1-2} to Eq. (\ref{EQ2}) we have the subcases below:\\
\underline{Subcase 1.1}: $f=\alpha(\chi_{1}-\chi_{2})$, $g_{1}=\dfrac{\chi_{1}+\chi_{2}}{2}$, where $\chi_{1},\chi_{2}\in\mathcal{M}(S)$  such that $\chi_{1}\neq\chi_{2}$ and $\alpha\in\mathbb{C}\setminus\{0\}$ is a constant. By a small computation, we get from (\ref{EQ5}) that $$\Psi_{m}(\delta\,a_{1}-2a+\delta\,a^{2})=-2\alpha\,\lambda_{2}(\chi_{1}-\chi_{2}).$$ As $\delta\,a_{1}-2a\in\mathcal{A}(S\setminus I_{m})$, we deduce, according to Proposition \ref{prop8-1}(\ref{(4)}), that $\delta\,a^{2}=0$. Since $a\neq0$ we get that $\delta=0$. So that the identity above reduces to $\Psi_{m}(a)=\alpha\,\lambda_{2}(\chi_{1}-\chi_{2})$, which implies, by applying Proposition \ref{prop8-1}(\ref{(2)})(\ref{(1)}), that $a=0$, contradicting the constraint on $a$. Hence, the system $(\ref{EQ1})-(\ref{EQ2})$ has no solution in this subcase.\\
\underline{Subcase 1.2}: $f=\Psi_{\chi}(A)$, $g_{1}=\chi$, where $0\neq\chi\in\mathcal{M}(S)$ and $0\neq A\in\mathcal{A}(S\setminus I_{\chi})$. Then, from  (\ref{EQ5}) we obtain $$\Psi_{m}(-\delta\,a_{1}+2\,a-\delta\,a^{2})=\Psi_{\chi}(2\lambda_{2}\,A).$$
Hence, according to Proposition \ref{prop8-1}(\ref{(6)}), we get that $\delta\,a=0$. As $a\neq0$ we deduce that $\delta=0$. Then, the identity above reduces to
\begin{equation}\label{EQ6}\Psi_{m}(a)=\Psi_{\chi}(\lambda_{2}\,A).\end{equation} As $m,\,\chi\in\mathcal{M}(S)\setminus\{0\}$ and $a\in\mathcal{A}(S\setminus I_{m})\setminus\{0\}$, we get, according to Proposition \ref{prop8-1-1}(\ref{(3)}), that $m=\chi$ and $a=\lambda_{2}\,A$.\\
So, taking (\ref{EQ3}) and (\ref{EQ4}) into account, we obtain $$f=\Psi_{m}(A),\,\,g_{1}=m$$
$$h=\dfrac{1}{2}\Psi_{m}(a_{1}+\lambda_{2}^{2}A^{2}),\,\,g_{2}=m.$$
The solution occurs in part (1).
\\\underline{Case 2}: $H=c^{2}(\mu-m)-c\Psi_{m}(a),\,\,\,\,G=m,\,\,\,\,\lambda_{2}\,F=c(\mu-m)$, where $c\in\mathbb{C}\setminus\{0\}$ is a constant, $m\neq0,\mu\in\mathcal{M}(S)$ such that $\mu\neq m$, and $0\neq a\in\mathcal{A}(S\setminus I_{m})$.\\
As in case 1, we have following subcases:\\
\underline{Subcase 2.1}: $f=\alpha\,(\chi_{1}-\chi_{2})$, $g_{1}=\dfrac{\chi_{1}+\chi_{2}}{2}$, where $\chi_{1},\chi_{2}\in\mathcal{M}(S)$  such that $\chi_{1}\neq\chi_{2}$ and $\alpha\in\mathbb{C}\setminus\{0\}$ is a constant. Then, we get from (\ref{EQ5}) that
 $$\Psi_{m}(\delta\,c\,a)=\alpha\,\lambda_{2}\,(\chi_{1}-\chi_{2})+(\delta\,c^{2}-c)(\mu-m).$$ Hence, by applying Proposition\ref{prop8-1}(\ref{(2)})(\ref{(1)}), we deduce that $\delta\,c\,a=0$. As $c\neq0$ and $a\neq0$ we get that $\delta=0$. So that the identity above reduces to
\begin{equation}\label{EQ7}\alpha\,\lambda_{2}\,\chi_{1}-\alpha\,\lambda_{2}\,\chi_{2}-c\,\mu+c\,m=0.\end{equation}
 Since $c\neq0$ we conclude, by \cite[Theorem 3.18]{Stetk1}, that the multiplicative functions $\chi_{1},\chi_{2},\mu$ and $m$ are not different. As $\chi_{1}\neq\chi_{2}$ and $\mu\neq m$ we have the following possibilities:\\
\underline{Subcase 2.1.1}: $\chi_{1}=\mu$. Then (\ref{EQ7}) implies $(\alpha\,\lambda_{2}-c)\,\chi_{1}-\alpha\,\lambda_{2}\,\chi_{2}+c\,m=0$. From which we get, according to \cite[Theorem 3.18]{Stetk1}, that $\chi_{2}=m$ because $c\neq0$ and $\chi_{1}\neq\chi_{2}$. Then $(\alpha\,\lambda_{2}-c)\,(\chi_{1}-\chi_{2})=0$. So, $c=\alpha\,\lambda_{2}$. Hence, taking (\ref{EQ3}) and (\ref{EQ4}) into account, we deduce that
 $$f=\alpha\,(\mu-m),\,\,g_{1}=\dfrac{\mu+m}{2}$$
$$h=\alpha^{2}\,\lambda_{2}^{2}(\mu-m)-\alpha\,\lambda_{2}\Psi_{m}(a),\,\,g_{2}=m.$$
The solution occurs in part (2).
\\\underline{Subcase 2.1.2}: $\chi_{1}=m$. Then (\ref{EQ7}) reduces to $(\alpha\,\lambda_{2}+c)\,\chi_{1}-\alpha\,\lambda_{2}\,\chi_{2}+c\,m=0$.
Hence, as in Subcase 2.1.1, we get that $\chi_{2}=m$ and $c=-\alpha\,\lambda_{2}$. By writing $-\alpha$ and $-a$ instead of $\alpha$ and $a$ respectively we go back to Subcase 2.1.1, and using that $\Psi_{m}$ is linear we get the same formulas for $(f,g_{1},h,g_{2})$ as the ones in Subcase 2.1.1.\\
\underline{Subcase 2.2}: $f=\Psi_{\chi}(A)$, $g_{1}=\chi$, where $0\neq\chi\in\mathcal{M}(S)$ and $0\neq A\in\mathcal{A}(S\setminus I_{\chi})$. Then, we get from (\ref{EQ5}) that
$$\Psi_{\chi}(\lambda_{2}\,A)+\Psi_{m}(-\delta\,c\,a)=(c-\delta\,c^{2})\mu-(c-\delta\,c^{2})m.$$
Since $\lambda_{2}\,A\in\mathcal{A}(S\setminus I_{\chi})$ and $-\delta\,c\,a\in\mathcal{A}(S\setminus I_{m})$ we deduce from the identity above, applying Proposition \ref{prop8-1}(3), that $\Psi_{\chi}(\lambda_{2}\,A)+\Psi_{m}(-\delta\,c\,a)=0$ and then $(c-\delta\,c^{2})(\mu-m)=0$. Since $\mu\neq m$ and $c\neq0$ we get that $\delta\,c=1$. So, taking (\ref{EQ3}) and (\ref{EQ4}) into account, we deduce the followings formulas for $(f,g_{1},h,g_{2})$:\\
$$f=\dfrac{1}{\lambda_{2}}\,\Psi_{m}(a),\,\,g_{1}=m$$
$$h=c^{2}(\mu-m)-c\Psi_{m}(a),\,\,g_{2}=\dfrac{\mu-m}{2}-\dfrac{1}{2c}\,\Psi_{m}(a).$$
The solution occurs in part (3).
\\\underline{Case 3}: $H=c(\mu-m)+c\,d\,\Psi_{m}(a),\,\,\,\,G=\dfrac{\mu+m}{2}-\dfrac{1}{2}d\,\Psi_{m}(a),\,\,\,\,\lambda_{2}\,F=\Psi_{m}(a)$, where $\mu,m\in\mathcal{M}(S)$ such that $m\neq\mu$ and $m\neq0$, $0\neq a\in\mathcal{A}(S\setminus I_{m})$ and $c,d\in\mathbb{C}\setminus\{0\}$ are constants such that $1-cd^{2}=0$.\\For the pair $(f,g_{1})$ we have the following subcases:\\
\underline{Subcase 3.1}: $f=\alpha\,(\chi_{1}-\chi_{2})$, $g_{1}=\dfrac{\chi_{1}+\chi_{2}}{2}$, where $\chi_{1},\chi_{2}\in\mathcal{M}(S)$  such that $\chi_{1}\neq\chi_{2}$ and $\alpha\in\mathbb{C}\setminus\{0\}$ is a constant. Then by a small computation , using that $\Psi_{m}$ is linear, (\ref{EQ5}) implies
$$\Psi_{m}((1-\delta\,c\,d)a)=\alpha\,\lambda_{2}(\chi_{1}-\chi_{2})+\delta\,c(\mu-m).$$  Hence, by applying Proposition \ref{prop8-1}(\ref{(2)})(\ref{(1)}), we deduce from the identity above that $(1-\delta\,c\,d)a=0$. As $a\neq0$ and $c\,d^{2}=1$ we get that $\delta=d$. then the identity above reduces to
\begin{equation}\label{EQ8}\alpha\,\lambda_{2}\,\chi_{1}-\alpha\,\lambda_{2}\,\chi_{2}+c\,d\,\mu-c\,d\,m=0.\end{equation}
Since $c\,d\neq0$, $\chi_{1}\neq\chi_{2}$ and $\mu\neq m$, we deduce by \cite[Theorem 3.18]{Stetk1} that we have the following subcases:\\
\underline{Subcase 3.1.1}: $\chi_{1}=\mu$. Then (\ref{EQ8}) reduces to $(\alpha\,\lambda_{2}+c\,d)\mu-\alpha\,\lambda_{2}\,\chi_{2}-c\,d\,m=0$, which implies, according to \cite[Theorem 3.18]{Stetk1}, that $\chi_{2}=m$ because $c\,d\neq0$ and $\chi_{1}\neq\chi_{2}$. So that $(\alpha\,\lambda_{2}+c\,d)(\mu-m)=0$. As $m\neq\mu$ we get that $c\,d=-\alpha\,\lambda_{2}$ which implies $\alpha=-\dfrac{1}{d\,\lambda_{2}}$ because $c\,d^{2}=1$.\\Hence, $$f=-\dfrac{1}{d\,\lambda_{2}}(\mu-m).$$
Taking (\ref{EQ3}), (\ref{EQ4}) and the formulas of $H$ and $G$ into account and using that $c\,d^{2}=1$, we derive
$$h=c(\mu-m)+c\,d\,\Psi_{m}(a)=c(\mu-m)+\dfrac{1}{d}\,\Psi_{m}(a)$$
and $$g_{2}=-\dfrac{c\,d^{2}}{2}(\mu-m)-\dfrac{c\,d^{3}}{2}\,\Psi_{m}(a)+\dfrac{\mu+m}{2}-\dfrac{1}{2}d\,\Psi_{m}(a)+d\,\Psi_{m}(a)=m.$$
The solution occurs in part (4).
\\\underline{Subcase 3.1.2}: $\chi_{1}=m$. Then (\ref{EQ8}) becomes $(\alpha\,\lambda_{2}-c\,d)m-\alpha\,\lambda_{2}\,\chi_{2}+c\,d\,\mu=0$. Since $c\,d\neq0$ and $\chi_{1}\neq\chi_{2}$, we get from the last identity, according to \cite[Theorem 3.18]{Stetk1}, that $\chi_{2}=\mu$. So that $(\alpha\,\lambda_{2}-c\,d)(m-\mu)=0$. As $m\neq\mu$ we get that $c\,d=\alpha\,\lambda_{2}$. By a small computation we obtain the same formulas for $(f,g_{1},h,g_{2})$ as the ones in Subcase 3.1.1.\\
\underline{Subcase 3.2}: $f=\Psi_{\chi}(A)$, $g_{1}=\chi$, where $0\neq\chi\in\mathcal{M}(S)$ and $0\neq A\in\mathcal{A}(S\setminus I_{\chi})$. Then, using that $\Psi_{m}$ is linear, we get from (\ref{EQ5}) that $$\Psi_{\chi}(\lambda_{2}\,A)+\Psi_{m}((\delta\,c\,d-1))a)=-\delta\,c(\mu-m).$$
So, according to Proposition \ref{prop8-1}(\ref{(3)}), we deduce from the identity above that $\Psi_{\chi}(\lambda_{2}\,A)$$\\$$+\Psi_{m}((\delta\,c\,d-1))a)=\delta\,c(\mu-m)=0$, which implies, seeing that $\mu\neq m$ and $c\neq0$, $$\Psi_{\chi}(\lambda_{2}\,A)=\Psi_{m}((1-\delta\,c\,d))a)\quad\text{and}\quad \delta=0.$$
So, $\Psi_{\chi}(\lambda_{2}\,A)=\Psi_{m}(a)$. As in Subcase 1.2 we get that $m=\chi$ and $a=\lambda_{2}\,A$.
So, in view of (\ref{EQ3}), (\ref{EQ4}) and $c\,d^{2}=1$, we obtain $$f=\dfrac{1}{\lambda_{2}}\Psi_{m}(a),\,\,g_{1}=m$$
$$h=c(\mu-m)+\dfrac{1}{d}\,\Psi_{m}(a),\,\,g_{2}=\dfrac{\mu+m}{2}-\dfrac{1}{2}d\,\Psi_{m}(a).$$
The solution occurs in part (5).
\\\underline{Case 4}: $$H=c\,\beta\,\chi_{1}+c\,(2-\beta)\,\chi_{2}-2\,c\,\chi_{3},\,\,G=\dfrac{1}{4}\,\beta\,\chi_{1}+\dfrac{1}{4}\,(2-\beta)\,\chi_{2}+\dfrac{1}{2}\,\chi_{3},\,\,\lambda_{2}\,F=\dfrac{\chi_{1}-\chi_{2}}{2\lambda},$$
where $\chi_{1},\chi_{2},\chi_{3}$ are different multiplicative functions on $S$ and $\lambda,\beta,c\in\mathbb{C}\setminus\{0\}$ are constants with $2\,c\,\lambda^{2}\,\beta\,(2-\beta)=1$.\\According to Proposition \ref{prop8-1-2} we have the following subcases for the pair $\{f,g_{1}\}$:\\
\underline{Subcase 4.1}: $f=\alpha\,(\mu_{1}-\mu_{2})$, $g_{1}=\dfrac{\mu_{1}+\mu_{2}}{2}$, where $\mu_{1},\mu_{2}\in\mathcal{M}(S)$  such that $\mu_{1}\neq\mu_{2}$ and $\alpha\in\mathbb{C}\setminus\{0\}$ is a constant. Then (\ref{EQ5}) implies
\begin{equation}\label{EQ9}\alpha\,\lambda_{2}(\mu_{1}-\mu_{2})+(\delta\,c\,\beta-\dfrac{1}{2\lambda})\chi_{1}+(\delta\,c(2-\beta)+\dfrac{1}{2\lambda})\chi_{2}-2\,\delta\,c\,\chi_{3}=0.\end{equation}
Since $\alpha\,\lambda_{2}\neq0$, $\mu_{1}\neq\mu_{2}$ and $\chi_{1},\chi_{2},\chi_{3}$ are different we get from (\ref{EQ9}), by applying \cite[Theorem 3.18]{Stetk1} two times successively, that we have the following subcases:\\
\underline{Subcase 4.1.1}: $\mu_{1}=\chi_{1}$ and $\mu_{2}=\chi_{2}$. Then (\ref{EQ9}) reduces to
\begin{equation}\label{EQ10}(\alpha\,\lambda_{2}+\delta\,c\,\beta-\dfrac{1}{2\lambda})\chi_{1}+(\delta\,c(2-\beta)+\dfrac{1}{2\lambda}-\alpha\,\lambda_{2})\chi_{2}-2\,\delta\,c\,\chi_{3}=0.\end{equation}
As $\chi_{1},\chi_{2},\chi_{3}$ are different multiplicative functions on $S$ then, according to \cite[Theorem 3.18]{Stetk1} and seeing that $c\neq0$, we get that $\delta=0$ and $\alpha=\dfrac{1}{2\,\lambda\,\lambda_{2}}$. Hence $$f=\dfrac{\chi_{1}-\chi_{2}}{2\,\lambda\,\lambda_{2}},\,\,g_{1}=\dfrac{\chi_{1}+\chi_{2}}{2}.$$
Moreover, in view of (\ref{EQ3}) and (\ref{EQ4}), we obtain
$$h=c\,\beta\,\chi_{1}+c\,(2-\beta)\,\chi_{2}-2\,c\,\chi_{3},\,\,g_{2}=\dfrac{\beta\,\chi_{1}+(2-\beta)\,\chi_{2}+2\,\chi_{3}}{4}.$$
The solution occurs in part (6).
\\\underline{Subcase 4.1.2}: $\mu_{1}=\chi_{1}$ and $\mu_{2}=\chi_{3}$. Then, from (\ref{EQ9}), we get\\
\begin{equation}\label{EQ11}(\alpha\,\lambda_{2}+\delta\,c\,\beta-\dfrac{1}{2\lambda})\chi_{1}+(\delta\,c(2-\beta)+\dfrac{1}{2\lambda})\chi_{2}-(\alpha\,\lambda_{2}+2\,\delta\,c)\chi_{3}=0.\end{equation}
Since the multiplicative functions $\chi_{1},\chi_{2},\chi_{3}$ are different we deduce, according to \cite[Theorem 3.18]{Stetk1}, that
\begin{equation}\label{EQ12}\alpha\,\lambda_{2}+\delta\,c\,\beta-\dfrac{1}{2\lambda}=0,\end{equation}
\begin{equation}\label{EQ13}\delta\,c(2-\beta)+\dfrac{1}{2\lambda}=0\end{equation} and
\begin{equation}\label{EQ14}\alpha\,\lambda_{2}+2\,\delta\,c=0.\end{equation}
The identity (\ref{EQ13}) implies $-2\,\lambda\,\delta\,c(2-\beta)=1$. As $2\,c\,\lambda^{2}\,\beta(2-\beta)=1$, we infer that $\delta=-\lambda\,\beta$. Substituting this (\ref{EQ14}) (or in (\ref{EQ12})) we get that $\alpha=\dfrac{1}{\lambda\,\lambda_{2}(2-\beta)}$. Hence, $f=\dfrac{1}{\lambda\,\lambda_{2}(2-\beta)}(\chi_{1}-\chi_{3})$ and $g_{1}=\dfrac{\chi_{1}+\chi_{3}}{2}$.\\
On the other hand, in view of (\ref{EQ3}) and (\ref{EQ4}), we get that
\begin{equation*}h=c\,\beta\,\chi_{1}+c\,(2-\beta)\,\chi_{2}-2\,c\,\chi_{3},\end{equation*}
and
\begin{equation*}\begin{split}&g_{2}=-\dfrac{1}{2}\delta^{2}[c\,\beta\,\chi_{1}+c\,(2-\beta)\,\chi_{2}-2\,c\,\chi_{3}]\\
&\quad\quad\,+\dfrac{1}{4}\,\beta\,\chi_{1}+\dfrac{1}{4}(2-\beta)\chi_{2}+\dfrac{1}{2}\,\chi_{3}-\dfrac{1}{2}\,\beta(\chi_{1}-\chi_{2})\\
&\quad=-\beta\,\dfrac{1+2\,c\,\lambda^{2}\,\beta^{2}}{4}\,\chi_{1}+((\beta-2)\,\dfrac{1+2\,c\,\lambda^{2}\,\beta^{2}}{4}+1)\,\chi_{2}+\dfrac{1+2\,c\,\lambda^{2}\,\beta^{2}}{2}\,\chi_{3}.\end{split}\end{equation*}
Since $2\,c\,\lambda^{2}\,\beta\,(2-\beta)=1$ we get from the identity above, by a small computation, that $$g_{2}=\dfrac{-\beta\,\chi_{1}+(2-\beta)\,\chi_{2}+2\,\chi_{3}}{2(2-\beta)}.$$
The solution occurs in part (7).
\\\underline{Subcase 4.1.3}: $\mu_{1}=\chi_{2}$ and $\mu_{2}=\chi_{1}$. Then (\ref{EQ9}) becomes
$$(-\alpha\,\lambda_{2}+\delta\,c\,\beta-\dfrac{1}{2\,\lambda})\chi_{1}+(\alpha\,\lambda_{2}+\delta\,c(2-\beta)+\dfrac{1}{2\,\lambda})\chi_{2}-2\,\delta\,c\,\chi_{3}=0.$$
As in subcase 4.1.1 we get that $\delta=0$
and $\alpha=-\dfrac{1}{2\,\lambda\,\lambda_{2}}$. So that we obtain the same solution as the one in Subcase 4.1.1.
\\\underline{Subcase 4.1.4}: $\mu_{1}=\chi_{2}$ and $\mu_{2}=\chi_{3}$. Then (\ref{EQ9}) reduces to
$$(\delta\,c\,\beta-\dfrac{1}{2\,\lambda})\chi_{1}+(\alpha\,\lambda_{2}+\delta\,c(2-\beta)+\dfrac{1}{2\,\lambda})\chi_{2}-(\alpha\,\lambda_{2}+2\,\delta\,c)\chi_{3}=0.$$
Since the multiplicative functions $\chi_{1},\chi_{2},\chi_{3}$ are different we deduce, according to \cite[Theorem 3.18]{Stetk1}, that $$\delta\,c\,\beta-\dfrac{1}{2\,\lambda}=\alpha\,\lambda_{2}+\delta\,c(2-\beta)+\dfrac{1}{2\,\lambda}=\alpha\,\lambda_{2}+2\,\delta\,c=0.$$ Taking the identity $2\,c\,\lambda^{2}\,\beta(2-\beta)=1$ into account, we get that $\delta=\lambda(2-\beta)$ and $\alpha=-\dfrac{1}{\lambda\,\lambda_{2}\,\beta}$. By similar computations to the ones in Subcase 4.1.2 we obtain
$$f=-\dfrac{1}{\lambda\,\lambda_{2}\,\beta}(\chi_{2}-\chi_{3}),\,\,g_{1}=\dfrac{\chi_{2}+\chi_{3}}{2},$$
$$h=c\,\beta\,\chi_{1}+c\,(2-\beta)\,\chi_{2}-2\,c\,\chi_{3},\,\,g_{2}=\dfrac{\beta\,\chi_{1}-(2-\beta)\,\chi_{2}+2\,\chi_{3}}{2\,\beta}.$$
The solution occurs in part (8).
\\\underline{Subcase 4.1.5}: $\mu_{1}=\chi_{3}$ and $\mu_{2}=\chi_{1}$. Then (\ref{EQ9}) implies that
$$(-\alpha\,\lambda_{2}+\delta\,c\,\beta-\dfrac{1}{2\lambda})\chi_{1}+(\delta\,c(2-\beta)+\dfrac{1}{2\lambda})\chi_{2}-(\alpha\,\lambda_{2}-2\,\delta\,c)\chi_{3}=0.$$
 Hence, by writing $-\alpha$ instead of $\alpha$, we go back to subcase 4.1.2 and we obtain the solution of the form (7) in Theorem \ref{thm1}. \\\underline{Subcase 4.1.6}: $\mu_{1}=\chi_{3}$ and $\mu_{2}=\chi_{2}$. Then (\ref{EQ9}) reduces to
$$(\delta\,c\,\beta-\dfrac{1}{2\lambda})\chi_{1}+(-\alpha\,\lambda_{2}+\delta\,c(2-\beta)+\dfrac{1}{2\lambda})\chi_{2}+(\alpha\,\lambda_{2}-2\,\delta\,c)\chi_{3}=0.$$
Hence, by writing $-\alpha$ instead of $\alpha$, we go back to subcase 4.1.4 and we obtain the solution of the form (8) in Theorem \ref{thm1}.
\\\underline{Subcase 4.2}: $f=\Psi_{\chi}(A)$, $g_{1}=\chi$, where $0\neq\chi\in\mathcal{M}(S)$ and $0\neq A\in\mathcal{A}(S\setminus I_{\chi})$. Using that $\Psi_{\chi}:\mathcal{F}(S\setminus I_{\chi},\mathbb{C})\to\mathcal{F}(S,\mathbb{C})$ is linear we get from (\ref{EQ5}) that $$\Psi_{\chi}(\lambda_{2}\,A)=(-\delta\,c\,\beta+\dfrac{1}{2\,\lambda})\chi_{1}+(-\delta\,c(2-\beta)-\dfrac{1}{2\,\lambda})\chi_{2}+2\delta\,c\,\chi_{3}.$$
So, according to Proposition \ref{prop8-1}(2)(1) and seeing that $\lambda_{2}\neq0$, we get that $A=0$ which contradicts the assumption on $A$. Hence, the system $(\ref{eq01})-(\ref{eq02})$ has no solution in this subcase.
\par Conversely if $f,g_{1},h$ and $g_{2}$ are of the
forms (1)-(8) in Theorem \ref{thm1} we check by elementary computations
that $f,g_{1},h$ and $g_{2}$ satisfy the system $(\ref{EQ1})-(\ref{EQ2})$, and
$f$ and $h$ are linearly independent. This completes the proof of Theorem \ref{thm1}.
\end{proof}
\section{Solution of the system $(\ref{eq01})-(\ref{eq02})$ when $\lambda_{1}\neq0$ and $\lambda_{2}\neq0$}
In this section we solve the following system, which we denote $(\ref{eq01})-(\ref{eq02})$:
for four unknown complex valued functions $f,g_{1},g_{2},h$ on a semigroup $S$ generated by its squares such that $f$ and $h$ are linearly independent, and $\lambda_{1},\lambda_{2}\in\mathbb{C}\setminus\{0\}$ are given constants.\\
The solutions of the system $(\ref{eq01})-(\ref{eq02})$ are given in Theorem \ref{thm2}.
\begin{thm}\label{thm2} The complex valued solutions of system $(\ref{eq01})-(\ref{eq02})$ such that  $f$ and $h$ are linearly independent are the following quadruplets $(f,g_{1},h,g_{2})$.
 \begin{enumerate}
\item\label{(t1)} $$f=\dfrac{1}{2}\,\Psi_{\chi}(A_{1}+A^{2}),$$ $$g_{1}=-\dfrac{1}{4}\delta\,\Psi_{\chi}(\delta\,A_{1}-4A)-\dfrac{1}{4}\delta^{2}\,\Psi_{\chi}(A^{2})+\chi,$$
    $$h=-\dfrac{1}{2\lambda_{1}}\Psi_{\chi}(\delta\,A_{1}-2\,A)-\dfrac{\delta}{2\,\lambda_{1}}\,\Psi_{\chi}(A^{2}),$$$$g_{2}=-\dfrac{1}{4}\delta\Psi_{\chi}(\delta\,A_{1}+2\,A)-\dfrac{1}{4}\delta^{2}\,\Psi_{\chi}(A^{2})+\chi,$$
    where $\chi\in\mathcal{M}(S)$ such that $\chi\neq0$, $A,A_{1}\in\mathcal{A}(S\setminus I_{\chi})$ such that $A\neq0$, and $\delta\in\mathbb{C}\setminus\{0\}$ is a constant.\\
\item\label{(t2)}$$f=\dfrac{\lambda_{1}}{\lambda_{1}\,\delta_{1}+\delta_{2}^{2}}(\dfrac{\lambda_{1}}{\lambda_{1}\,\delta_{1}+\delta_{2}^{2}}(\mu-\chi)-\Psi_{\chi}(A),$$
$$g_{1}=\dfrac{\lambda_{1}\,\delta_{1}}{2(\lambda_{1}\,\delta_{1}+\delta_{2}^{2})}((1+\dfrac{\delta_{2}^{2}}{\lambda_{1}\,\delta_{1}+\delta_{2}^{2}})(\mu-\chi)+\Psi_{\chi}(A))+\chi,$$
$$h=(\dfrac{\lambda_{2}}{\lambda_{2}\,\delta_{2}+\delta_{1}^{2}})^{2}(\mu-\chi),$$ $$g_{2}=\dfrac{\lambda_{2}\,\delta_{2}}{2(\lambda_{2}\,\delta_{2}+\delta_{1}^{2})}((1+\dfrac{\delta_{1}^{2}}{\lambda_{2}\,\delta_{2}+\delta_{1}^{2}})(\mu-\chi)-\delta_{1}\,\Psi_{\chi}(A))+\chi,$$
where $\chi\neq0,\mu\in\mathcal{M}(S)$ and $A\in\mathcal{A}(S\setminus I_{\chi})$ such that $\mu\neq\chi$ and $A\neq0$, and $\delta_{1},\delta_{2}\in\mathbb{C}\setminus\{0\}$ are constants such that $\lambda_{2}\,\delta_{2}+\delta_{1}^{2}\neq0$ and $\lambda_{1}\,\delta_{1}+\delta_{2}^{2}\neq0$.\\
\item $$f=c^{2}(\mu-\chi)-c\,\Psi_{\chi}(A),$$
$$g_{1}=\dfrac{(1-(\delta\,c-1)^{2})\mu+(1+(\delta\,c-1)^{2})\chi}{2}+\dfrac{1}{2}\,\delta^{2}\,c\,\Psi_{\chi}(A,$$
$$h=\dfrac{\lambda_{2}^{2}\,c^{2}}{\delta^{2}}(\mu-\chi)+\dfrac{\delta\,c}{\lambda_{1}}\,\Psi_{\chi}(A),$$
$$g_{2}=\dfrac{(1-\delta^{2}\,c^{2})\mu+(1+\delta^{2}\,c^{2})\chi}{2}-\dfrac{\delta^{3}(c-\delta)^{2}}{2\,c\,\lambda_{1}\,\lambda_{2}^{2}}\,\Psi_{\chi}(A),$$
where $\chi\neq0,\mu\in\mathcal{M}(S)$ and $A\in\mathcal{A}(S\setminus I_{\chi})$ such that $\mu\neq\chi$ and $A\neq0$, and $\delta,c\in\mathbb{C}\setminus\{0\}$ are constants such that $c\,\delta^{3}-\delta^{2}+c\,\lambda_{1}\,\lambda_{2}^{2}=0$.
\item\label{(t4)}$$f=c(\mu-\chi)+c\,d\,\Psi_{\chi}(A),$$ $$g_{1}=\dfrac{(d^{2}-\delta^{2})\mu+(d^{2}+\delta^{2})\chi}{2\,d^{2}}+\dfrac{d^{2}+\delta^{2}}{2\,d}\Psi_{\chi}(A),$$
    $$h=-\dfrac{\delta}{d^{2}\,\lambda_{1}}(\mu-\chi)+\dfrac{d-\delta}{d\,\lambda_{1}}\Psi_{\chi}(A),$$
    $$g_{2}=\dfrac{(\delta\,d^{2}+\lambda_{1}\,\lambda_{2}^{2})\mu+(\delta\,d^{2}-\lambda_{1}\,\lambda_{2}^{2})\chi}{2\,\delta\,d^{2}}-\dfrac{\lambda_{1}\,\lambda_{2}^{2}}{2\,d(d-\delta)}\,\Psi_{\chi},$$
where  $\chi\neq0,\mu\in\mathcal{M}(S)$ and $A\in\mathcal{A}(S\setminus I_{\chi})$ such that $\mu\neq\chi$ and $A\neq0$, and $c,d,\delta\in\mathbb{C}\setminus\{0\}$ are constants such that $d\neq\delta$, $1-c\,d^{2}=0$.
\item\label{(t5)} There exist different multiplicative functions $\chi_{1},\chi_{2},\chi_{3}$ on $S$ and constants $\alpha,\beta,\gamma,\lambda,\delta_{1},\delta_{2},c,d\in\mathbb{C}\setminus\{0\}$ satisfying $2\,c\,\lambda^{2}\,d\,(2-d)=1$ and $2\,\alpha\,\gamma^{2}\,\beta\,(2-\beta)=1$ such that
$$f=F_{1},\, g_{1}=-\dfrac{1}{2}\delta_{1}^{2}\,F_{1}+G_{1}+\delta_{1}\lambda_{1}\,H_{1}$$
$$h=H_{2},\,g_{2}=-\dfrac{1}{2}\delta_{2}^{2}\,H_{2}+G_{2}+\delta_{2}\lambda_{2}\,F_{2},$$
where
$$F_{1}=c\,d\,\chi_{1}+c\,(2-d)\,\chi_{2}-2\,c\,\chi_{3},\,\,G_{1}=\dfrac{1}{4}\,d\,\chi_{1}+\dfrac{1}{4}\,(2-d)\,\chi_{2}+\dfrac{1}{2}\,\chi_{3},$$$$\lambda_{1}\,H_{1}=\dfrac{\chi_{1}-\chi_{2}}{2\,\lambda},$$
and the triplet $(H_{2},G_{2},\lambda_{2}\,F_{2})$ is of the following forms:
\\(i) $$H_{2}=\alpha\,\beta\,\chi_{1}+\alpha\,(2-\beta)\,\chi_{2}-2\,\alpha\,\chi_{3},\,\,G_{2}=\dfrac{1}{4}\,\beta\,\chi_{1}+\dfrac{1}{4}\,(2-\beta)\,\chi_{2}+\dfrac{1}{2}\,\chi_{3},$$$$\lambda_{2}\,F_{2}=\dfrac{\chi_{1}-\chi_{2}}{2\,\gamma},$$
 with $$\delta_{1}=\dfrac{\gamma\,\lambda_{2}}{\lambda},\,\delta_{2}=\dfrac{\lambda\,\lambda_{1}}{\gamma}, \, \alpha=-\dfrac{c\,\gamma\,\lambda_{2}}{\lambda\,\lambda_{1}},\,\beta=d-\dfrac{1}{2\,c\,\gamma\,\lambda_{2}}.$$
\\(ii)
$$H_{2}=\alpha\,\beta\,\chi_{2}+\alpha\,(2-\beta)\,\chi_{1}-2\,\alpha\,\chi_{3},\,\,G_{2}=\dfrac{1}{4}\,\beta\,\chi_{2}+\dfrac{1}{4}\,(2-\beta)\,\chi_{1}+\dfrac{1}{2}\,\chi_{3},$$$$\lambda_{2}\,F_{2}=\dfrac{\chi_{1}-\chi_{2}}{2\,\gamma},$$
with $$\delta_{1}=-\dfrac{\gamma\,\lambda_{2}}{\lambda},\,\delta_{2}=-\dfrac{\lambda\,\lambda_{1}}{\gamma},\,\alpha=\dfrac{c\,\gamma\,\lambda_{2}}{\lambda\,\lambda_{1}},\,\beta=2-d-\dfrac{1}{2\,c\,\gamma\,\lambda_{2}}.$$
\\(iii)
$$H_{2}=\alpha\,\beta\,\chi_{1}+\alpha\,(2-\beta)\,\chi_{3}-2\,\alpha\,\chi_{2},\,\,G_{2}=\dfrac{1}{4}\,\beta\,\chi_{1}+\dfrac{1}{4}\,(2-\beta)\,\chi_{3}+\dfrac{1}{2}\,\chi_{2},$$$$\lambda_{2}\,F_{2}=\dfrac{\chi_{1}-\chi_{3}}{2\,\gamma},$$ with $$\delta_{1}=\dfrac{4\,c\,\gamma\,\lambda_{2}-1}{2\,c(2-d)\,\lambda},\,\delta_{2}=\dfrac{(2-d)\,\lambda\,\lambda_{1}}{2\,\gamma},\, \alpha=\dfrac{c\,\gamma\,\lambda_{2}}{\lambda\,\lambda_{1}},\,\beta=\dfrac{1-2\,c\,d\,\gamma\,\lambda_{2}}{c(2-d)\gamma\,\lambda_{2}}.$$
\\(iv)
$$H_{2}=\alpha\,\beta\,\chi_{3}+\alpha\,(2-\beta)\,\chi_{1}-2\,\alpha\,\chi_{2},\,\,G_{2}=\dfrac{1}{4}\,\beta\,\chi_{3}+\dfrac{1}{4}\,(2-\beta)\,\chi_{1}+\dfrac{1}{2}\,\chi_{2},$$$$\lambda_{2}\,F_{2}=\dfrac{\chi_{3}-\chi_{1}}{2\,\gamma},$$
 with $$\delta_{1}=-\dfrac{1+4\,c\,\gamma\,\lambda_{2}}{2\,c(2-d)\,\lambda},\,\delta_{2}=-\dfrac{(2-d)\,\lambda\,\lambda_{1}}{2\,\gamma},\,\alpha=-\dfrac{c\,\gamma\,\lambda_{2}}{\lambda\,\lambda_{1}}, ,\,\beta=\dfrac{1+4\,c\,\gamma\,\lambda_{2}}{c(2-d)\gamma\,\lambda_{2}}.$$
\\(v)
$$H_{2}=\alpha\,\beta\,\chi_{2}+\alpha\,(2-\beta)\,\chi_{3}-2\,\alpha\,\chi_{1},\,\,G_{2}=\dfrac{1}{4}\,\beta\,\chi_{2}+\dfrac{1}{4}\,(2-\beta)\,\chi_{3}+\dfrac{1}{2}\,\chi_{1},$$$$\lambda_{2}\,F_{2}=\dfrac{\chi_{2}-\chi_{3}}{2\,\gamma},$$
 with $$\delta_{1}=\dfrac{1-4\,c\,\gamma\,\lambda_{2}}{2\,c\,d\,\lambda},\,\delta_{2}=-\dfrac{d\,\lambda\,\lambda_{1}}{2\,\gamma},\, \alpha=-\dfrac{c\,\gamma\,\lambda_{2}}{\lambda\,\lambda_{1}},\,\beta=\dfrac{1-2\,c(2-d)\gamma\,\lambda_{2}}{c\,d\,\gamma\,\lambda_{2}}.$$
\\(vi)
$$H_{2}=\alpha\,\beta\,\chi_{3}+\alpha\,(2-\beta)\,\chi_{2}-2\,\alpha\,\chi_{1},\,\,G_{2}=\dfrac{1}{4}\,\beta\,\chi_{3}+\dfrac{1}{4}\,(2-\beta)\,\chi_{2}+\dfrac{1}{2}\,\chi_{1},$$$$\lambda_{2}\,F_{2}=\dfrac{\chi_{3}-\chi_{2}}{2\,\gamma},$$
with $$\delta_{1}=\dfrac{1+4\,c\,\gamma\,\lambda_{2}}{2\,c\,d\,\lambda},\,\delta_{2}=\dfrac{d\,\lambda\,\lambda_{1}}{2\,\gamma},\, \alpha=\dfrac{c\,\gamma\,\lambda_{2}}{\lambda\,\lambda_{1}},\,\beta=\dfrac{1+4\,c\,\gamma\,\lambda_{2}}{c\,d\,\gamma\,\lambda_{2}}.$$
\end{enumerate}
\end{thm}
\begin{proof} Applying Proposition \ref{prop8-2} we get from functional equations (\ref{eq01}) and (\ref{eq02}) that
\begin{equation}\label{eq3} f=F_{1},\end{equation}
\begin{equation}\label{eq4} g_{1}=-\dfrac{1}{2}\delta_{1}^{2}\,F_{1}+G_{1}+\delta_{1}\lambda_{1}\,H_{1},\end{equation}
\begin{equation}\label{eq5} \lambda_{1}\,h=-\delta_{1}\,F_{1}+\lambda_{1}\,H_{1},\end{equation}
\begin{equation}\label{eq6} h=H_{2},\end{equation}
\begin{equation}\label{eq7} g_{2}=-\dfrac{1}{2}\delta_{2}^{2}\,H_{2}+G_{2}+\delta_{2}\lambda_{2}\,F_{2},\end{equation}
\begin{equation}\label{eq8} \lambda_{2}\,f=-\delta_{2}\,H_{2}+\lambda_{2}\,F_{2},\end{equation}
where $\delta_{1},\delta_{2}\in\mathbb{C}$ are constants, and the triplets $(F_{1},G_{1},\lambda_{1}\,H_{1})$ and $(H_{2},G_{2},\lambda_{2}\,F_{2})$ are of the forms in Proposition \ref{prop8-2}(\ref{(8)}),(\ref{(9)}),(\ref{(10)}),(\ref{(11)}) with the same constraints.\\
Hence, from (\ref{eq5}) and (\ref{eq6}), and (\ref{eq3}) and (\ref{eq8})  we get respectively
\begin{equation}\label{eq9} H_{1}-H_{2}=\alpha_{1}\,F_{1}\end{equation}
and
\begin{equation}\label{eq10} F_{1}-F_{2}=\alpha_{2}\,H_{2}\end{equation}
where $\alpha_{1}:=\dfrac{\delta_{1}}{\lambda_{1}}$ and $\alpha_{2}:=-\dfrac{\delta_{2}}{\lambda_{2}}$.\\
According to Proposition \ref{prop8-2}, we split the discussion into the following cases:\\
\underline{Case 1}: $F_{1}=\dfrac{1}{2}\Psi_{\chi}(A_{1}+A^{2}),\,\,\,\,G_{1}=\chi,\,\,\,\,\lambda_{1}\,H_{1}=\Psi_{\chi}(A)$, where $0\neq \chi\in\mathcal{M}(S)$, $A\neq0,A_{1}\in\mathcal{A}(S\setminus I_{\chi})$.\\
For the triplet $(H_{2},G_{2},\lambda_{2}\,F_{2})$ we have the following possibilities:\\
\underline{Subcase 1.1}: $H_{2}=\dfrac{1}{2}\Psi_{m}(a_{1}+a^{2}),\,\,\,\,G_{2}=m,\,\,\,\,\lambda_{2}\,F_{2}=\Psi_{m}(a)$, where $0\neq m\in\mathcal{M}(S)$, $a\neq0,a_{1}\in\mathcal{A}(S\setminus I_{m})$. Then, we get from (\ref{eq9}) and (\ref{eq10}) that
\begin{equation}\label{eq11}\dfrac{1}{\lambda_{1}}\Psi_{\chi}(A)-\dfrac{1}{2}\Psi_{m}(a_{1}+a^{2})=\dfrac{\alpha_{1}}{2}\Psi_{\chi}(A_{1}+A^{2})\end{equation}
and
\begin{equation}\label{eq12}\dfrac{1}{2}\Psi_{\chi}(A_{1}+A^{2})-\dfrac{1}{\lambda_{2}}\Psi_{m}(a)=\dfrac{\alpha_{2}}{2}\Psi_{m}(a_{1}+a^{2}).\end{equation}
Now, multiplying (\ref{eq12}) by $\alpha_{1}$, and adding the identity obtained to (\ref{eq11}), and using the linearity of $\Psi_{\chi}$ and $\Psi_{m}$, we get that
\begin{equation}\label{eq13}\Psi_{\chi}(\dfrac{2}{\lambda_{1}}\,A)=\Psi_{m}(\dfrac{2\,\alpha_{1}}{\lambda_{2}}\,a+(1+\alpha_{1}\,\alpha_{2})a_{1}+(1+\alpha_{1}\,\alpha_{2})a^{2}).\end{equation}
So that, according to Proposition  \ref{prop8-1}(\ref{(6)}) and taking into account that $a\neq0$, we deduce that $1+\alpha_{1}\,\alpha_{2}=0$. Then $\delta_{1}\,\delta_{2}=\lambda_{1}\,\lambda_{2}$. Moreover, (\ref{eq13}) reduces to $\Psi_{\chi}(\dfrac{1}{\lambda_{1}}\,A)=\Psi_{m}(\dfrac{\alpha_{1}}{\lambda_{2}}\,a)=\Psi_{m}(\dfrac{\delta_{1}}{\lambda_{1}\,\lambda_{2}}\,a)=\Psi_{m}(\dfrac{1}{\delta_{2}}\,a)$. So, according to Proposition \ref{prop8-1-1}(3) we derive
\begin{equation}\label{eq14}m=\chi\end{equation} and
\begin{equation}\label{eq15}a=\dfrac{\lambda_{2}}{\delta_{1}}\,A=\dfrac{\delta_{2}}{\lambda_{1}}\,A.\end{equation} Substituting this in (\ref{eq11}) and using Proposition  \ref{prop8-2}(\ref{(1)}) we get by a small computation, that $\dfrac{2}{\lambda_{1}}\,A-a_{1}-\dfrac{\delta_{1}}{\lambda_{1}}\,A_{1}=\dfrac{\delta_{1}^{3}+\lambda_{1}\lambda_{2}^{2}}{\lambda_{1}\delta_{1}^{2}}\,A^{2}$. So, the function $\dfrac{\delta_{1}^{3}+\lambda_{1}\lambda_{2}^{2}}{\lambda_{1}\delta_{1}^{2}}\,A^{2}$ is additive on $S\setminus I_{\chi}$. As $A\neq0$ we get that $\delta_{1}^{3}+\lambda_{1}\lambda_{2}^{2}=0$. Then $\delta_{1}^{2}+\lambda_{2}\delta_{2}=0$, $\delta_{2}^{2}+\lambda_{1}\delta_{1}=0$ ( because $\delta_{1}\,\delta_{2}=\lambda_{1}\,\lambda_{2}$ ) and
\begin{equation}\label{eq16}a_{1}=\dfrac{2}{\lambda_{1}}\,A-\dfrac{\delta_{1}}{\lambda_{1}}\,A_{1}.\end{equation}
So, using the identities (\ref{eq3}), (\ref{eq4}), (\ref{eq6}), (\ref{eq7}), (\ref{eq15}) and (\ref{eq16}), and taking the formulas of $F_{1},G_{1},F_{2},G_{2},H_{1},H_{2}$ in this subcase into account, we obtain
\begin{equation*}f=\dfrac{1}{2}\,\Psi_{\chi}(A_{1}+A^{2}),\end{equation*}
\begin{equation*}g_{1}=-\dfrac{1}{4}\delta_{1}\,\Psi_{\chi}(\delta_{1}\,A_{1}-4A)-\dfrac{1}{4}\delta_{1}^{2}\,\Psi_{\chi}(A^{2})+\chi,\end{equation*}
\begin{equation*}\begin{split}&h=\dfrac{1}{2}\Psi_{\chi}(a_{1}+a^{2})\\
&\,\,\,\,=\dfrac{1}{2}\Psi_{\chi}(\dfrac{2}{\lambda_{1}}\,A-\dfrac{\delta_{1}}{\lambda_{1}}\,A_{1}+\dfrac{\delta_{2}^{2}}{\lambda_{1}^{2}}\,A^{2})\\
&\,\,\,\,=-\dfrac{1}{2\lambda_{1}}\Psi_{\chi}(\delta_{1}\,A_{1}-2\,A)-\dfrac{\delta_{1}}{2\,\lambda_{1}}\,\Psi_{\chi}(A^{2}),\end{split}\end{equation*}
and
\begin{equation*}\begin{split}&g_{2}=-\dfrac{1}{4}\delta_{2}^{2}\Psi_{\chi}(a_{1}+a^{2})+\chi+\delta_{2}\Psi_{\chi}(a)\\
&\,\,\,\,\,\,=-\dfrac{1}{4}\delta_{2}^{2}\Psi_{\chi}(\dfrac{2}{\lambda_{1}}\,A-\dfrac{\delta_{1}}{\lambda_{1}}\,A_{1}+\dfrac{\lambda_{2}^{2}}{\delta_{1}^{2}}\,A^{2})+\chi+\Psi_{\chi}(\dfrac{\delta_{2}^{2}}{\lambda_{1}}\,A)\\
&\,\,\,\,\,\,=-\dfrac{1}{4}\delta_{1}\Psi_{\chi}(\delta_{1}\,A_{1}+2\,A)-\dfrac{1}{4}\delta_{1}^{2}\,\Psi_{\chi}(A^{2})+\chi.\end{split}\end{equation*}
The solution occurs in part (1).
\\\underline{Subcase 1.2}: $H_{2}=c^{2}(\mu-m)-c\Psi_{m}(a),\,\,\,\,G_{2}=m,\,\,\,\,\lambda_{2}\,F_{2}=c(\mu-m)$, where $c\in\mathbb{C}\setminus\{0\}$ is a constant, $m\neq0,\mu\in\mathcal{M}(S)$ such that $\mu\neq m$, and $0\neq a\in\mathcal{A}(S\setminus I_{m})$. Then, (\ref{eq10}) implies
$$\Psi_{\chi}(\dfrac{1}{2}\,A_{1}+\dfrac{1}{2}\,A^{2})+\Psi_{m}(\alpha_{2}\,c\,a)=(\dfrac{c}{\lambda_{2}}+\alpha_{2}\,c^{2})(\mu-m).$$
So, according to Proposition \ref{prop8-1}(5) we get that $A=0$, which contradicts the assumption on $A$. Hence, the system $(\ref{eq01})-(\ref{eq02})$ has no solution in subcase.\\
\underline{Subcase 1.3}: $H_{2}=c(\mu-m)+c\,d\,\Psi_{m}(a),\,\,\,\,G_{2}=\dfrac{\mu+m}{2}-\dfrac{1}{2}d\,\Psi_{m}(a),\,\,\,\,\lambda_{2}\,F_{2}=\Psi_{m}(a)$, where $\mu,m\in\mathcal{M}(S)$ such that $m\neq\mu$ and $m\neq0$, $0\neq a\in\mathcal{A}(S\setminus I_{m})$ and $c,d\in\mathbb{C}\setminus\{0\}$ are constants such that $1-cd^{2}=0$. Then, by a small computation, we get from (\ref{eq10}) that $\Psi_{\chi}(\dfrac{1}{2}\,A_{1}+\dfrac{1}{2}\,A^{2})-\Psi_{m}((\dfrac{1}{\lambda_{2}}+\alpha_{2}^{2}\,cd)a)=\alpha_{2}^{2}\,c(\mu-m)$. So, by applying Proposition \ref{prop8-1}(5) we get that $A=0$, which contradicts the assumption on $A$. Hence, the system $(\ref{eq01})-(\ref{eq02})$ has no solution in subcase.\\
\underline{Subcase 1.4}:
$$H_{2}=c\,\beta\,\chi_{1}+c\,(2-\beta)\,\chi_{2}-2\,c\,\chi_{3},\,\,G_{2}=\dfrac{1}{4}\,\beta\,\chi_{1}+\dfrac{1}{4}\,(2-\beta)\,\chi_{2}+\dfrac{1}{2}\,\chi_{3},\,\,\lambda_{2}\,F_{2}=\dfrac{\chi_{1}-\chi_{2}}{2\,\lambda},$$
where $\chi_{1},\chi_{2},\chi_{3}$ are different multiplicative functions on $S$ and $\lambda,\beta,c\in\mathbb{C}\setminus\{0\}$ are constants with $2\,c\,\lambda^{2}\,\beta\,(2-\beta)=1$. Then, (\ref{eq10}) becomes
$$\Psi_{\chi}(\dfrac{1}{2}\,A_{1}+\dfrac{1}{2}\,A^{2})=\alpha_{2}^{2}\,c\,\beta\,\chi_{1}+\alpha_{2}^{2}\,c\,(2-\beta)\,\chi_{2}-2\,\alpha_{2}^{2}\,\,c\,\chi_{3}+\dfrac{\chi_{1}-\chi_{2}}{2\,\lambda\,\lambda_{2}}.$$
So, according to Proposition \ref{prop8-1}(4), we obtain $A=0$, contradicting the constraint on $A$. It follows that the system $(\ref{eq01})-(\ref{eq02})$ has no solution in subcase.\\
\underline{Case 2}: $F_{1}=c^{2}(\mu-\chi)-c\Psi_{\chi}(A),\,\,\,\,G_{1}=\chi,\,\,\,\,\lambda_{1}\,H_{1}=c(\mu-\chi)$, where $c\in\mathbb{C}\setminus\{0\}$ is a constant, $\chi\neq0,\mu\in\mathcal{M}(S)$ with $\mu\neq\chi$, and $0\neq A\in\mathcal{A}(S\setminus I_{\chi})$.\\
For the triplet $(H_{2},G_{2},\lambda_{2}\,F_{2})$ we have the following possibilities:\\
\underline{Subcase 2.1}: $H_{2}=\dfrac{1}{2}\Psi_{m}(a_{1}+a^{2}),\,\,\,\,G_{2}=m,\,\,\,\,\lambda_{2}\,F_{2}=\Psi_{m}(a)$, where $0\neq m\in\mathcal{M}(S)$, $a\neq0,a_{1}\in\mathcal{A}(S\setminus I_{m})$.\\By interchanging equations (\ref{eq01}) and (\ref{eq02}) we go back to Subcase 1.2. Then, the system $(\ref{eq01})-(\ref{eq02})$ has no solution in this subcase.\\
\underline{Subcase 2.2}: $H_{2}=d^{2}(m_{1}-m_{2})-d\,\Psi_{m_{2}}(a),\,\,\,\,G_{2}=m_{2},\,\,\,\,\lambda_{2}\,F_{2}=d(m_{1}-m_{2})$, where $d\in\mathbb{C}\setminus\{0\}$ is a constant, $m_{2}\neq0,m_{1}\in\mathcal{M}(S)$ with $m_{1}\neq m_{2}$, and $0\neq a\in\mathcal{A}(S\setminus I_{m_{2}})$. Then, from (\ref{eq9}) and (\ref{eq10}) we derive respectively
\begin{equation}\label{eq17}\Psi_{m_{2}}(d\,a)+\Psi_{\chi}(\alpha_{1}\,c\,A)=d^{2}(m_{1}-m_{2})+(\alpha_{1}\,c^{2}-\dfrac{c}{\lambda_{1}})\mu-(\alpha_{1}\,c^{2}-\dfrac{c}{\lambda_{1}})\chi.\end{equation}
\begin{equation}\label{eq18}\Psi_{\chi}(c\,A)-\Psi_{m_{2}}(\alpha_{2}\,d\,a)=c^{2}(\mu-\chi)-(\alpha_{2}\,d^{2}+\dfrac{d}{\lambda_{2}})m_{1}+(\alpha_{2}\,d^{2}+\dfrac{d}{\lambda_{2}})m_{2}.\end{equation}
Since $c\neq0$, $d\neq0$, $a\neq0$ and $A\neq0$ we deduce from (\ref{eq17}) and (\ref{eq18}), by applying Propositions \ref{prop8-1}(3)(1) and \ref{prop8-1-1}(3), that
 \begin{equation}\label{eq018}m_{2}=\chi,\end{equation}
$c\,A=\alpha_{2}\,d\,a$ and $d\,a=-\alpha_{1}\,c\,A$, then $A=-\alpha_{1}\,\alpha_{2}\,A$. Hence, $\alpha_{1}\,\alpha_{2}=-1$. So,
\begin{equation}\label{eq019}a=-\dfrac{\delta_{1}\,c}{\lambda_{1}\,d}A\end{equation}
and
\begin{equation}\label{eq19} \delta_{1}\,\delta_{2}=\lambda_{1}\,\lambda_{2}.\end{equation} Moreover, (\ref{eq17}) and (\ref{eq18}) reduce to
$$(\dfrac{c}{\lambda_{1}}-\alpha_{1}\,c^{2}-d^{2})\chi+(\alpha_{1}\,c^{2}-\dfrac{c}{\lambda_{1}})\mu+d^{2}\,m_{1}=0$$
and
$$(\dfrac{d}{\lambda_{2}}+\alpha_{2}\,d^{2}-c^{2})\chi+c^{2}\,\mu-(\alpha_{2}\,d^{2}+\dfrac{d}{\lambda_{2}})m_{1}=0.$$
Since $\chi\neq\mu,\,\chi\neq m_{1}$, $c\neq0$ and $d\neq0$ we derive from the identities above, according to \cite[Theorem 3.18]{Stetk1}, that
\begin{equation}\label{eq0019}\mu=m_{1},\end{equation}
\begin{equation}\label{eq20}\dfrac{c}{\lambda_{1}}-\alpha_{1}\,c^{2}-d^{2}=0\end{equation}
and
\begin{equation}\label{eq21}\dfrac{d}{\lambda_{2}}+\alpha_{2}\,d^{2}-c^{2}=0.\end{equation}
Now, we multiply (\ref{eq20}) by $\alpha_{2}$ by using that $\alpha_{1}\,\alpha_{2}=-1$. Then, by adding the identity obtained to (\ref{eq21}) and taking into account that $\delta_{1}\,\delta_{2}=\lambda_{1}\,\lambda_{2}\neq0$, we derive that
\begin{equation}\label{eq22}d=\dfrac{\delta_{2}}{\lambda_{1}}c=\dfrac{\lambda_{2}}{\delta_{1}}c.\end{equation}
Substituting this in (\ref{eq019}) and (\ref{eq20}) we get respectively, by a small computation, that
\begin{equation}\label{eq022}a=-\dfrac{\delta_{1}}{\delta_{2}}A\end{equation}
and $(\delta_{1}\,\lambda_{1}+\delta_{2}^{2})c=\lambda_{1}$. So, $\delta_{1}\,\lambda_{1}+\delta_{2}^{2}\neq0$ because $\lambda_{1}\neq0$. Then,
\begin{equation}\label{eq23}c=\dfrac{\delta_{1}^{2}}{\delta_{1}^{3}+\lambda_{1}\,\lambda_{2}^{2}}.\end{equation}
Hence, taking (\ref{eq19}) and (\ref{eq22}) into account, we obtain
\begin{equation}\label{eq24}d=\dfrac{\lambda_{2}\,\delta_{1}}{\delta_{1}^{3}+\lambda_{1}\,\lambda_{2}^{2}}.\end{equation}
So, using the identities (\ref{eq3}), (\ref{eq4}), (\ref{eq6}), (\ref{eq7}), (\ref{eq018}), (\ref{eq0019}), (\ref{eq022}), (\ref{eq23}) and (\ref{eq24}), and taking the formulas of $F_{1},G_{1},F_{2},G_{2},H_{1},H_{2}$ in this subcase into account, we obtain
\begin{equation*}f=\dfrac{\lambda_{1}}{\lambda_{1}\,\delta_{1}+\delta_{2}^{2}}(\dfrac{\lambda_{1}}{\lambda_{1}\,\delta_{1}+\delta_{2}^{2}}(\mu-\chi)-\Psi_{\chi}(A),\end{equation*}
\begin{equation*}\begin{split}&g_{1}=-\dfrac{1}{2}\delta_{1}^{2}\dfrac{\lambda_{1}^{2}}{(\lambda_{1}\,\delta_{1}+\delta_{2}^{2})^{2}}(\mu-\chi)+\dfrac{1}{2}\delta_{1}^{2}\dfrac{\lambda_{1}}{\lambda_{1}\,\delta_{1}+\delta_{2}^{2}}\Psi_{\chi}(A)+\dfrac{\lambda_{1}\,\delta_{1}}{\lambda_{1}\,\delta_{1}+\delta_{2}^{2}}(\mu-\chi)+\chi\\
&\,\,\,\,\,\,=\dfrac{\lambda_{1}\,\delta_{1}}{2(\lambda_{1}\,\delta_{1}+\delta_{2}^{2})}((1+\dfrac{\delta_{2}^{2}}{\lambda_{1}\,\delta_{1}+\delta_{2}^{2}})(\mu-\chi)+\Psi_{\chi}(A))+\chi,\end{split}\end{equation*}
\begin{equation*}h=(\dfrac{\lambda_{2}}{\lambda_{2}\,\delta_{2}+\delta_{1}^{2}})^{2}(\mu-\chi)\end{equation*}
and
\begin{equation*}\begin{split}&g_{2}=-\dfrac{1}{2}\delta_{2}^{2}\dfrac{\lambda_{2}^{2}}{(\lambda_{2}\,\delta_{2}+\delta_{1}^{2})^{2}}(\mu-\chi)+\dfrac{1}{2}\delta_{2}^{2}\dfrac{\lambda_{2}}{\lambda_{2}\,\delta_{2}+\delta_{1}^{2}}\Psi_{\chi}(a)+\dfrac{\lambda_{2}\,\delta_{2}}{\lambda_{2}\,\delta_{2}+\delta_{1}^{2}}(\mu-\chi)+\chi\\
&\,\,\,\,\,\,=\dfrac{\lambda_{2}\,\delta_{2}}{2(\lambda_{2}\,\delta_{2}+\delta_{1}^{2})}((1+\dfrac{\delta_{1}^{2}}{\lambda_{2}\,\delta_{2}+\delta_{1}^{2}})(\mu-\chi)-\delta_{1}\,\Psi_{\chi}(A))+\chi,\end{split}\end{equation*}
The solution occurs in part (2).
\\\underline{Subase 2.3}: $H_{2}=\alpha(m_{1}-m_{2})+\alpha\,\beta\,\Psi_{m_{2}}(a),\,\,\,\,G_{2}=\dfrac{m_{1}+m_{2}}{2}-\dfrac{1}{2}\beta\,\Psi_{m_{2}}(a),\,\,\,\,\lambda_{2}\,F_{2}=\Psi_{m_{2}}(a)$, where $m_{1},m_{2}\in\mathcal{M}(S)$ such that $m_{1}\neq m_{2}$ and $m_{2}\neq0$, $0\neq a\in\mathcal{A}(S\setminus I_{m_{2}})$ and $\alpha,\beta\in\mathbb{C}\setminus\{0\}$ are constants such that $1-\alpha\,\beta^{2}=0$. Then, from (\ref{eq9}) and (\ref{eq10}), we get that
\begin{equation}\label{eq25}\Psi_{\chi}(\alpha_{1}\,cA)-\Psi_{m_{2}}(\alpha\,\beta\,a)=(\alpha_{1}\,c^{2}-\dfrac{c}{\lambda_{1}})\mu-(\alpha_{1}\,c^{2}-\dfrac{c}{\lambda_{1}})\chi+\alpha(m_{1}-m_{2})\end{equation}
and
\begin{equation}\label{eq26}\Psi_{m_{2}}(\alpha_{2}\,\alpha\,\beta+\dfrac{1}{\lambda_{2}})a)+\Psi_{\chi}(cA)=c^{2}(\mu-\chi)-\alpha_{2}\,\alpha(m_{1}-m_{2}).\end{equation}
Seeing that $c,\alpha,\beta\in\mathbb{C}\setminus\{0\}$, $a\neq0$ and $A\neq0$, and using Propositions  \ref{prop8-1}(3)(1) and \ref{prop8-1-1}(3), we deduce from (\ref{eq25}) and (\ref{eq26}) that
\begin{equation}\label{eq27}m_{2}=\chi,\end{equation}
\begin{equation}\label{eq28}m_{1}=\mu,\end{equation}
\begin{equation}\label{eq29}\alpha\,\beta\,a=\alpha_{1}\,cA,\end{equation}
\begin{equation}\label{eq30}(\alpha_{2}\,\alpha\,\beta+\dfrac{1}{\lambda_{2}})a=-c\,A,\end{equation}
\begin{equation}\label{eq31}\alpha_{1}\,c^{2}-\dfrac{c}{\lambda_{1}}+\alpha=0\end{equation}
and
\begin{equation}\label{eq32}c^{2}=\alpha_{2}\,\alpha.\end{equation}
Multiplying (\ref{eq29}) by $\alpha_{2}$ and substituting the identity obtained in (\ref{eq30}) derive $a=-c\,\lambda_{2}(\alpha_{1}\,\alpha_{2}+1)A$. Then $$a=c\,\dfrac{\delta_{1}\,\delta_{2}-\lambda_{1}\,\lambda_{2}}{\lambda_{1}}\,A.$$
Multiplying (\ref{eq29}) by $\beta$ and using that $\alpha\,\beta^{2}=1$ we obtain
\begin{equation}\label{eq32-1} a=\dfrac{c\,\beta\,\delta_{1}}{\lambda_{1}}\,A.\end{equation}
As $A\neq0$ we deduce from identities above that
\begin{equation}\label{eq32-2} \beta=\dfrac{\delta_{1}\,\delta_{2}-\lambda_{1}\,\lambda_{2}}{\delta_{1}}.\end{equation}
On the other hand, we get from (\ref{eq31}) and (\ref{eq32}) that $\alpha_{1}\,\alpha_{2}\,\alpha-\dfrac{c}{\lambda_{1}}+\alpha=0$. Then $\alpha(1+\alpha_{1}\,\alpha_{2})=\dfrac{c}{\lambda_{1}}$, which implies that $\alpha(\lambda_{1}\,\lambda_{2}-\delta_{1}\,\delta_{2})=c\,\lambda_{2}$. So, since $c\,\lambda_{2}\neq0$, we deduce that $\delta_{1}\,\delta_{2}\neq\lambda_{1}\,\lambda_{2}$ and $\alpha=\dfrac{c\,\lambda_{2}}{\lambda_{1}\,\lambda_{2}-\delta_{1}\,\delta_{2}}$. Then, using (\ref{eq32-2}) we get that
\begin{equation}\label{eq32-3} \alpha\,\beta=-\dfrac{c\,\lambda_{2}}{\delta_{1}}.\end{equation}
As $\alpha\,\beta^{2}=1$ we drive from (\ref{eq32-3}) that
\begin{equation}\label{eq32-4} \beta=-\dfrac{\delta_{1}}{c\,\lambda_{2}}\end{equation}
and
\begin{equation}\label{eq32-5} \alpha=\dfrac{c^{2}\,\lambda_{2}^{2}}{\delta_{1}^{2}}.\end{equation}
Hence, from (\ref{eq32}) and (\ref{eq32-5}) we deduce that $c^{2}=-c^{2}\,\dfrac{\delta_{2}\,\lambda_{2}}{\delta_{1}^{2}}$. As $c\neq0$, we get that
\begin{equation}\label{eq32-6} \delta_{2}=-\dfrac{\delta_{1}^{2}}{\lambda_{2}}.\end{equation}
Substituting (\ref{eq32-4}) in (\ref{eq32-1}) we obtain
\begin{equation}\label{eq32-7} a=-\dfrac{\delta_{1}^{2}}{\lambda_{1}\,\lambda_{2}}\,A.\end{equation}
By using the identities (\ref{eq3}), (\ref{eq4}), (\ref{eq7}) and the formulas of the triplet $(F_{1},G_{1},\lambda_{1}\,H_{1})$ in this case we derive
\begin{equation*} f=c^{2}(\mu-\chi)-c\,\Psi_{\chi}(A)\end{equation*}
and
\begin{equation*} \begin{split}&g_{1}=-\dfrac{1}{2}\,\delta_{1}^{2}\,c^{2}(\mu-\chi)+\dfrac{1}{2}\,\delta_{1}^{2}\,c\,\Psi_{\chi}(A)+\chi+\delta_{1}\,c(\mu-\chi)\\
&\,\,\,\,\,\,=\dfrac{(1-(\delta_{1}\,c-1)^{2})\mu+(1+(\delta_{1}\,c-1)^{2})\chi}{2}+\dfrac{1}{2}\,\delta_{1}^{2}\,c\,\Psi_{\chi}(A).\end{split}\end{equation*}
Now, from (\ref{eq6}), (\ref{eq27}), (\ref{eq28}), (\ref{eq32-3}), (\ref{eq32-5}) and (\ref{eq32-7}) we deduce that
\begin{equation*} h=\dfrac{\lambda_{2}^{2}\,c^{2}}{\delta_{1}^{2}}(\mu-\chi)+\dfrac{\delta_{1}\,c}{\lambda_{1}}\,\Psi_{\chi}(A).\end{equation*}
On the other hand, by using (\ref{eq7}), (\ref{eq27}), (\ref{eq28}), (\ref{eq32-4}), (\ref{eq32-6}) and (\ref{eq32-7}) we obtain
\begin{equation*} \begin{split}&g_{2}=-\dfrac{2\,\delta_{1}^{4}}{\lambda_{2}^{2}}(\dfrac{\lambda_{2}^{2}\,c^{2}}{\delta_{1}^{2}}(\mu-\chi)+\dfrac{\delta_{1}\,c}{\lambda_{1}}\,\Psi_{\chi}(A))+\dfrac{\mu+\chi}{2}-\dfrac{1}{2}\,\beta\,\Psi_{\chi}(a)+\delta_{2}\,\Psi_{\chi}(a)\\
&\,\,\,\,\,\,=-\dfrac{\delta_{1}^{2}\,c^{2}}{2}(\mu-\chi)+\dfrac{\mu+\chi}{2}-\dfrac{1}{2}(\dfrac{\delta_{1}^{5}\,c}{\lambda_{1}\,\lambda_{2}^{2}}+\dfrac{\delta_{1}^{3}}{c\,\lambda_{1}\,\lambda_{2}^{2}}-\dfrac{2\,\delta_{1}^{4}}{\lambda_{1}\,\lambda_{2}^{2}})\Psi_{\chi}(A),\end{split}\end{equation*}
which reduces to
$$g_{2}=\dfrac{(1-\delta_{1}^{2}\,c^{2})\mu+(1+\delta_{1}^{2}\,c^{2})\chi}{2}-\dfrac{\delta_{1}^{3}(c-\delta_{1})^{2}}{2\,c\,\lambda_{1}\,\lambda_{2}^{2}}\,\Psi_{\chi}(A).$$
Moreover, from (\ref{eq32-2}), (\ref{eq32-4}) and (\ref{eq32-6}) we get that $c$ and $\delta_{1}$ satisfy the equation $c\,\delta_{1}^{3}-\delta_{1}^{2}+c\,\lambda_{1}\,\lambda_{2}^{2}=0$.
The solution occurs in part (3).
\\\underline{Subase 2.4}:
$$H_{2}=\alpha\,\beta\,\chi_{1}+\alpha\,(2-\beta)\,\chi_{2}-2\,\alpha\,\chi_{3},\,\,G_{2}=\dfrac{1}{4}\,\beta\,\chi_{1}+\dfrac{1}{4}\,(2-\beta)\,\chi_{2}+\dfrac{1}{2}\,\chi_{3},$$$$\lambda_{2}\,F_{2}=\dfrac{\chi_{1}-\chi_{2}}{2\,\lambda},$$
where $\chi_{1},\chi_{2},\chi_{3}$ are different multiplicative functions on $S$ and $\lambda,\beta,\alpha\in\mathbb{C}\setminus\{0\}$ are constants with $2\,\alpha\,\lambda^{2}\,\beta\,(2-\beta)=1$. Then, by applying  Proposition \ref{prop8-1}(2)(1) and seeing that $c\neq0$,  we get from (\ref{eq10}) that $A=0$, which contradicts the assumption on $A$. Hence, the system $(\ref{eq01})-(\ref{eq02})$ has no solution in this subcase.\\
\underline{Case 3}: $F_{1}=c(\mu-\chi)+c\,d\,\Psi_{\chi}(A),\,\,\,\,G_{1}=\dfrac{\mu+\chi}{2}-\dfrac{1}{2}d\,\Psi_{\chi}(A),\,\,\,\,\lambda_{1}\,H_{1}=\Psi_{\chi}(A)$, where $\mu,\chi\in\mathcal{M}(S)$ such that $\chi\neq\mu$ and $\chi\neq0$, $0\neq A\in\mathcal{A}(S\setminus I_{\chi})$ and $c,d\in\mathbb{C}\setminus\{0\}$ are constants such that $1-cd^{2}=0$. Then, we have the following subcases:\\
\underline{Subcase 3.1}: $H_{2}=\dfrac{1}{2}\Psi_{m}(a_{1}+a^{2}),\,\,\,\,G_{2}=m,\,\,\,\,\lambda_{2}\,F_{2}=\Psi_{m}(a)$, where $0\neq m\in\mathcal{M}(S)$, $a\neq0,a_{1}\in\mathcal{A}(S\setminus I_{m})$. As in subcase 1-3, we deduce that the system $(\ref{eq01})-(\ref{eq02})$ has no solution in this subcase.\\
\underline{Subcase 3.2}: $H_{2}=\alpha^{2}(m_{1}-m_{2})-\alpha\,\Psi_{m_{2}}(a),\,\,\,\,G_{2}=m_{2},\,\,\,\,\lambda_{2}\,F_{2}=\alpha(m_{1}-m_{2})$, where $\alpha\in\mathbb{C}\setminus\{0\}$ is a constant, $m_{2}\neq0,m_{1}\in\mathcal{M}(S)$ with $m_{1}\neq m_{2}$, and $0\neq a\in\mathcal{A}(S\setminus I_{m_{2}})$. So, we go back to subcase 2.3.\\
\underline{Subcase 3.3}:
$$H_{2}=\alpha(m_{1}-m_{2})+\alpha\,\beta\,\Psi_{m_{2}}(a),\,\,\,\,G_{2}=\dfrac{m_{1}+m_{2}}{2}-\dfrac{1}{2}\beta\,\Psi_{m_{2}}(a),$$ $$\lambda_{2}\,F_{2}=\Psi_{m_{2}}(a),$$ where $m_{1},m_{2}\in\mathcal{M}(S)$ such that $m_{1}\neq m_{2}$ and $m_{2}\neq0$, $0\neq a\in\mathcal{A}(S\setminus I_{m_{2}})$ and $\alpha,\beta\in\mathbb{C}\setminus\{0\}$ are constants such that $1-\alpha\,\beta^{2}=0$.\\
$2\,\alpha\,\lambda^{2}\,\beta\,(2-\beta)=1$.
Using (\ref{eq9}) and (\ref{eq10}) we obtain
\begin{equation}\label{eq33}\Psi_{\chi}(c\,d\,A)-\Psi_{m_{2}}((\dfrac{1}{\lambda_{2}}+\alpha\,\beta\,\alpha_{2})a)=\alpha\,\alpha_{2}(m_{1}-m_{2})-c(\mu-\chi)\end{equation}
 and
\begin{equation}\label{eq34}\Psi_{\chi}((\dfrac{1}{\lambda_{1}}-c\,d\,\alpha_{1})A)-\Psi_{m_{2}}(\alpha\,\beta\,a)=c\,\alpha_{1}(\mu-\chi)+\alpha(m_{1}-m_{2}).\end{equation}
Recall that $c,d\in\mathbb{C}\setminus\{0\}$ and $A\neq0$. Then, according to Propositions \ref{prop8-1}(3)(1) and \ref{prop8-1-1}(3), we deduce from (\ref{eq33}) and (\ref{eq34}) that\\
\begin{equation}\label{eq35}\alpha\,\alpha_{2}(m_{1}-m_{2})=c(\mu-\chi),\end{equation}
\begin{equation}\label{eq36}c\,\alpha_{1}(\mu-\chi)=-\alpha(m_{1}-m_{2}),\end{equation}
\begin{equation}\label{eq37}\chi=m_{2},\end{equation}
\begin{equation}\label{eq38}c\,d\,A=(\dfrac{1}{\lambda_{2}}+\alpha\,\beta\,\alpha_{2})a\end{equation}
and
\begin{equation}\label{eq39}(\dfrac{1}{\lambda_{1}}-c\,d\,\alpha_{1})A=\alpha\,\beta\,a.\end{equation}
Multiplying (\ref{eq36}) by $\alpha_{2}$ and taking (\ref{eq35}) into account, we deduce that $c(\alpha_{1}\,\alpha_{2}+1)(\mu-\chi)=0$. Therefore $\alpha_{1}\,\alpha_{2}=-1$ because $c\neq0$ and $\mu\neq\chi$. Hence,
\begin{equation}\label{eq40}\delta_{1}\,\delta_{2}=\lambda_{1}\,\lambda_{2}.\end{equation}
Now, when we multiply (\ref{eq38}) by $d$ and (\ref{eq39}) by $\beta$ and using $c\,d^{2}=\alpha\,\beta^{2}=1$ we get that $(\dfrac{d}{\lambda_{2}}+\alpha\,\beta\,d\,\alpha_{2})(\dfrac{\beta}{\lambda_{1}}-c\,d\,\beta\,\alpha_{1})A=A$. So, we find after reduction that $\beta(d-\delta_{1})=d\,\delta_{2}$. As $d\,\delta_{2}\neq0$ we get that $d\neq\delta_{1}$. So, taking (\ref{eq40}) into account, we obtain
\begin{equation}\label{eq41}\beta=\dfrac{d\,\lambda_{1}\,\lambda_{2}}{\delta_{1}(d-\delta_{1})}.\end{equation}
On the other hand, in view of (\ref{eq37}), the identity (\ref{eq36}) can be written as follows $c\,\alpha_{1}\,\mu-(c\,\alpha_{1}+\alpha)\chi+\alpha\,m_{1}=0$. As $c\neq0$, $\mu,\chi\in\mathcal{M}(S)$ and $\mu\neq\chi$, we derive according to \cite[Theorem 3.18]{Stetk1}, that
\begin{equation}\label{eq42}\mu=m_{1}\end{equation}
and $c\,\alpha_{1}+\alpha=0$. Hence, using (\ref{eq40}) we get that
\begin{equation}\label{eq43}c=-\dfrac{\lambda_{1}}{\delta_{1}}\,\alpha=-\dfrac{\delta_{2}}{\lambda_{2}}\,\alpha.\end{equation}
Now, multiplying (\ref{eq39}) by $\beta$ and using (\ref{eq41}), $\alpha\,\beta^{2}=1$ and $\alpha_{1}=\dfrac{\delta_{1}}{\lambda_{1}}$ we deduce by a small computation
\begin{equation}\label{eq44}a=\dfrac{\delta_{2}}{\lambda_{1}}A=\dfrac{\lambda_{2}}{\delta_{1}}A.\end{equation}
Now, using (\ref{eq3}), (\ref{eq4}), (\ref{eq6}), (\ref{eq7}), formulas of $(F_{1},G_{1},\lambda_{1}\,H_{1})$ and $(H_{2},G_{2},\lambda_{2}\,H_{2})$, (\ref{eq37}) and (\ref{eq42}) an elementary computation shows that
\begin{equation*}f=c(\mu-\chi)+c\,d\,\Psi_{\chi}(A),\end{equation*}
\begin{equation*}\begin{split}&g_{1}=\dfrac{(1-c\,\delta_{1}^{2})\mu+(1+c\,\delta_{1}^{2})\chi}{2}+\dfrac{d(1+c\,\delta_{1}^{2})}{2}\Psi_{\chi}(A)\\
&\,\,\,\,\,=\dfrac{(d^{2}-\delta_{1}^{2})\mu+(d^{2}+\delta_{1}^{2})\chi}{2\,d^{2}}+\dfrac{d^{2}+\delta_{1}^{2}}{2\,d}\Psi_{\chi}(A),\end{split}\end{equation*}
\begin{equation*}\begin{split}&h=\alpha(\mu-\chi)+\alpha\,\beta\,\Psi_{\chi}(a)\\
&\,\,\,\,\,=-\dfrac{\delta_{1}}{\lambda_{1}}c(\mu-\chi)+\dfrac{\delta_{1}(d-\delta_{1})}{d\,\lambda_{1}\,\lambda_{2}}\times\dfrac{\lambda_{2}}{\delta_{1}}\Psi_{\chi}(A)\\
&\,\,\,\,\,=-\dfrac{\delta_{1}}{d^{2}\,\lambda_{1}}(\mu-\chi)+\dfrac{d-\delta_{1}}{d\,\lambda_{1}}\Psi_{\chi}(A)\end{split}\end{equation*}
and
\begin{equation*}\begin{split}&g_{2}=-\dfrac{\delta_{1}\,\,\delta_{2}^{2}}{2\,d^{2}\,\lambda_{1}}(\mu-\chi)-\dfrac{\delta_{2}^{2}(d-\delta_{1})}{2\,d\,\lambda_{1}}\Psi_{\chi}(A)+\dfrac{\mu+\chi}{2}-\dfrac{1}{2}\,\beta\,\Psi_{\chi}(a)+\delta_{2}\,\Psi_{\chi}(a)\\
&\,\,\,\,\,=\dfrac{(d^{2}+\delta_{2}\,\lambda_{2})\mu+(d^{2}-\delta_{2}\,\lambda_{2})\chi}{2\,d^{2}}-\dfrac{\lambda_{1}\,\lambda_{2}^{2}(d-\delta_{1})}{2\,d\,\delta_{1}^{2}}\,\Psi_{\chi}(A)-\dfrac{d\,\lambda_{1}\,\lambda_{2}^{2}}{\delta_{1}^{2}(d-\delta_{1})}\,\Psi_{\chi}(A)\\
&\,\,\,\,\,\,\,\,\,\,\,+\dfrac{\lambda_{1}\,\lambda_{2}^{2}}{\delta_{1}^{2}}\,\Psi_{\chi}(A),\end{split}\end{equation*}
from which we derive, by a small computation and taking (\ref{eq40}) into account, that
$$g_{2}=\dfrac{(\delta_{1}\,d^{2}+\lambda_{1}\,\lambda_{2}^{2})\mu+(\delta_{1}\,d^{2}-\lambda_{1}\,\lambda_{2}^{2})\chi}{2\,\delta_{1}\,d^{2}}-\dfrac{\lambda_{1}\,\lambda_{2}^{2}}{2\,d(d-\delta_{1})}\,\Psi_{\chi}.$$ The solution occurs in part (4).
\\\underline{Subcase 3.4}:
$$H_{2}=\alpha\,\beta\,m_{1}+\alpha\,(2-\beta)\,m_{2}-2\,\alpha\,m_{3},\,\,G_{2}=\dfrac{1}{4}\,\beta\,m_{1}+\dfrac{1}{4}\,(2-\beta)\,m_{2}+\dfrac{1}{2}\,m_{3},$$$$\lambda_{2}\,F_{2}=\dfrac{m_{1}-m_{2}}{2\,\lambda},$$
where $m_{1},m_{2},m_{3}$ are different multiplicative functions on $S$ and $\lambda,\beta,\alpha\in\mathbb{C}\setminus\{0\}$ are constants with $2\,\alpha\,\lambda^{2}\,\beta\,(2-\beta)=1$. Then, by applying  Proposition \ref{prop8-1}(2)(1) and seeing that $c\neq0$,  we get from (\ref{eq10}) that $A=0$, which contradicts the assumption on $A$. Hence, the system $(\ref{eq01})-(\ref{eq02})$ has no solution in this subcase.\\
\underline{Case 4}:
$$F_{1}=c\,d\,\chi_{1}+c\,(2-d)\,\chi_{2}-2\,c\,\chi_{3},\,\,G_{1}=\dfrac{1}{4}\,d\,\chi_{1}+\dfrac{1}{4}\,(2-d)\,\chi_{2}+\dfrac{1}{2}\,\chi_{3},$$$$\lambda_{1}\,H_{1}=\dfrac{\chi_{1}-\chi_{2}}{2\,\lambda},$$
where $\chi_{1},\chi_{2},\chi_{3}$ are different multiplicative functions on $S$ and $\lambda,c,d\in\mathbb{C}\setminus\{0\}$ are constants such that $2\,c\,\lambda^{2}\,d\,(2-d)=1$.\\
The following subcases: \\
\underline{Subcase 4.1}:
$H_{2}=\dfrac{1}{2}\Psi_{m}(a_{1}+a^{2}),\,\,\,\,G_{2}=m,\,\,\,\,\lambda_{2}\,F_{2}=\Psi_{m}(a)$, where $0\neq m\in\mathcal{M}(S)$, $a\neq0,a_{1}\in\mathcal{A}(S\setminus I_{m})$,\\
\underline{Subcase 4.2}:
$H_{2}=\alpha^{2}(m_{1}-m_{2})-\alpha\,\Psi_{m_{2}}(a),\,\,\,\,G_{2}=m_{2},\,\,\,\,\lambda_{2}\,F_{2}=\alpha(m_{1}-m_{2})$, where $\alpha\in\mathbb{C}\setminus\{0\}$ is a constant, $m_{2}\neq0,m_{1}\in\mathcal{M}(S)$ with $m_{1}\neq m_{2}$, and $0\neq a\in\mathcal{A}(S\setminus I_{m_{2}})$,\\
\underline{Subcase 4.3}:
$$H_{2}=\alpha(m_{1}-m_{2})+\alpha\,\beta\,\Psi_{m_{2}}(a),\,\,\,\,G_{2}=\dfrac{m_{1}+m_{2}}{2}-\dfrac{1}{2}\beta\,\Psi_{m_{2}}(a),$$ $$\lambda_{2}\,F_{2}=\Psi_{m_{2}}(a),$$ where $m_{1},m_{2}\in\mathcal{M}(S)$ such that $m_{1}\neq m_{2}$ and $m_{2}\neq0$, $0\neq a\in\mathcal{A}(S\setminus I_{m_{2}})$ and $\alpha,\beta\in\mathbb{C}\setminus\{0\}$ are constants such that $1-\alpha\,\beta^{2}=0$,\\
go back, respectively, to subcases 1.4, 2.4 and 3.4 in this section. Hence, the system $(\ref{eq01})-(\ref{eq02})$ has no solution in each of  subcases 4.1, 4.2 and 4.3.\\
Now, we deal with the following subcase:\\
\underline{Subcase 4.4}:
$$H_{2}=\alpha\,\beta\,m_{1}+\alpha\,(2-\beta)\,m_{2}-2\,\alpha\,m_{3},\,\,G_{2}=\dfrac{1}{4}\,\beta\,m_{1}+\dfrac{1}{4}\,(2-\beta)\,m_{2}+\dfrac{1}{2}\,m_{3},$$$$\lambda_{2}\,F_{2}=\dfrac{m_{1}-m_{2}}{2\,\gamma},$$
where $m_{1},m_{2},m_{3}$ are different multiplicative functions on $S$ and $\gamma,\beta,\alpha\in\mathbb{C}\setminus\{0\}$ are constants with $2\,\alpha\,\gamma^{2}\,\beta\,(2-\beta)=1$.\\
Using (\ref{eq9}) and (\ref{eq10}) we obtain
\begin{equation}\begin{split}&\label{eq45}(\dfrac{1}{2\,\lambda\,\lambda_{1}}-c\,d\,\alpha_{1})\chi_{1}-(\dfrac{1}{2\,\lambda\,\lambda_{1}}+c(2-d)\alpha_{1})\chi_{2}\\
&+2\,c\,\alpha_{1}\,\chi_{3}-\alpha\,\beta\,m_{1}-\alpha\,(2-\beta)\,m_{2}+2\,\alpha\,m_{3}=0.\end{split}\end{equation}
and
\begin{equation}\begin{split}&\label{eq46} c\,d\,\chi_{1}+c\,(2-d)\,\chi_{2}-2\,c\,\chi_{3}-(\dfrac{1}{2\,\gamma\,\lambda_{2}}+\alpha\,\beta\,\alpha_{2})m_{1}\\ &+(\dfrac{1}{2\,\gamma\,\lambda_{2}}-\alpha(2-\beta)\alpha_{2})m_{2}+2\,\alpha\,\alpha_{2}\,m_{3}=0.\end{split}\end{equation}
When we multiply (\ref{eq46}) by $\alpha_{1}$ and add the identity obtained to (\ref{eq45}) we derive, by a small computation
\begin{equation}\begin{split}&\label{eq47}\dfrac{1}{2\,\lambda\,\lambda_{1}}\,\chi_{1}-\dfrac{1}{2\,\lambda\,\lambda_{1}}\,\chi_{2}-(\dfrac{\alpha_{1}}{2\,\gamma\,\lambda_{2}}+\alpha\,\beta(1+\alpha_{1}\,\alpha_{2}))m_{1}\\
&+(\dfrac{\alpha_{1}}{2\,\gamma\,\lambda_{2}}-\alpha(2-\beta)(1+\alpha_{1}\,\alpha_{2}))m_{2}+2\,\alpha(1+\alpha_{1}\,\alpha_{2})m_{3}=0.\end{split}\end{equation}
Similarly, by multiplying  (\ref{eq45}) by $\alpha_{2}$ and subtracting  (\ref{eq46}) from the identity obtained we deduce
\begin{equation}\begin{split}&\label{eq48}(\dfrac{\alpha_{2}}{2\,\lambda\,\lambda_{1}}-c\,d(1+\alpha_{1}\,\alpha_{2}))\chi_{1}-(\dfrac{\alpha_{2}}{2\,\lambda\,\lambda_{1}}+c(2-d)(1+\alpha_{1}\,\alpha_{2}))\chi_{2}\\
&+2\,c(1+\alpha_{1}\,\alpha_{2})\chi_{3}+\dfrac{1}{2\,\gamma\,\lambda_{2}}\,m_{1}-\dfrac{1}{2\,\gamma\,\lambda_{2}}\,m_{2}=0.\end{split}\end{equation}
Notice that $\chi_{1},\chi_{2},\chi_{3}, m_{1}, m_{2},m_{3}\in\mathcal{M}(S)$, $\chi_{1},\chi_{2},\chi_{3}$ are different and $m_{1}, m_{2},m_{3}$ also, and $\alpha\neq0$. So, according to \cite[Theorem 3.18]{Stetk1}, in view of (\ref{eq45}) we have the following subcases:
\\\underline{Subcase 4.4.1}: $\chi_{1}=m_{1}$ and $\chi_{2}=m_{2}$. Then,  from (\ref{eq46}) we derive that $\chi_{3}=m_{3}$ because $c\neq0$. Moreover, (\ref{eq48}) implies that $1+\alpha_{1}\,\alpha_{2}=0$. Hence, $\delta_{1}\,\delta_{2}=\lambda_{1}\,\lambda_{2}$. So that (\ref{eq48}) becomes $(\dfrac{\alpha_{2}}{2\,\lambda\,\lambda_{1}}+\dfrac{1}{2\,\gamma\,\lambda_{2}})(\chi_{1}-\chi_{2})=0$. As $\chi_{1}\neq\chi_{2}$ we get that $\dfrac{\alpha_{2}}{\lambda\,\lambda_{1}}=-\dfrac{1}{\gamma\,\lambda_{2}}$, which implies that
\begin{equation}\label{eq49}\delta_{2}=\dfrac{\lambda\,\lambda_{1}}{\gamma}.\end{equation}
As $\delta_{1}\,\delta_{2}=\lambda_{1}\,\lambda_{2}$ we get from (\ref{eq49}) that
\begin{equation}\label{eq50}\delta_{1}=\dfrac{\gamma\,\lambda_{2}}{\lambda}.\end{equation}
Now, since $\chi_{1}=m_{1},\,\chi_{2}=m_{2},\,\chi_{3}=m_{3}$ are different multiplicative function on $S$, we deduce from (\ref{eq45}), by applying \cite[Theorem 3.18]{Stetk1}, that
\begin{equation}\label{eq51}\dfrac{1}{2\,\lambda\,\lambda_{1}}-c\,d\,\alpha_{1}-\alpha\,\beta=0,\end{equation}
\begin{equation}\label{eq52}\dfrac{1}{2\,\lambda\,\lambda_{1}}+c(2-d)\,\alpha_{1}+\alpha(2-\beta)=0\end{equation}
and
\begin{equation}\label{eq53}c\,\alpha_{1}+\alpha=0.\end{equation}
Hence, from (\ref{eq53}) we derive
\begin{equation}\label{eq54}\alpha=-\dfrac{c\,\gamma\,\lambda_{2}}{\lambda\,\lambda_{1}}\end{equation}
When we  substitute this in (\ref{eq51}) and seeing that $\alpha_{1}=\dfrac{\delta_{1}}{\lambda_{1}}=\dfrac{\gamma\,\lambda_{2}}{\lambda\,\lambda_{1}}$ we derive that $\beta=\dfrac{c\,d\,\gamma\,\lambda\,\lambda_{1}\,\lambda_{2}}{c\,\gamma\,\lambda\,\lambda_{1}\,\lambda_{2}}-\dfrac{\lambda\,\lambda_{1}}{2\,\lambda\,\lambda_{1}\,c\gamma\,\lambda_{2}}$.
Hence, $\beta=d-\dfrac{1}{2\,c\,\gamma\,\lambda_{2}}$.\\
The solution occurs in part (5)(i).
\\\underline{Subcase 4.4.2}: $\chi_{1}=m_{2}$ and $\chi_{2}=m_{1}$. Then,  from (\ref{eq46}) we derive that $\chi_{3}=m_{3}$ because $c\neq0$. Proceeding as in subcase 4.1.1 we obtain $\delta_{1}=-\dfrac{\gamma\,\lambda_{2}}{\lambda}$, $\delta_{2}=-\dfrac{\lambda\,\lambda_{1}}{\gamma}$, $\alpha=\dfrac{c\,\gamma\,\lambda_{2}}{\lambda\,\lambda_{1}}$ and $\beta=\dfrac{2\,c(2-d)\,\gamma\,\lambda_{2}-1}{2\,c\,\gamma\,\lambda_{2}}=2-d-\dfrac{1}{2\,c\,\gamma\,\lambda_{2}}$.\\
The solution occurs in part (5)(ii).
\\\underline{Subcase 4.4.3}: $\chi_{1}=m_{1}$ and $\chi_{2}=m_{3}$. Then, from (\ref{eq46}) we derive that $\chi_{3}=m_{2}$ because $c\neq0$. So, from (\ref{eq47}) we get that
\begin{equation}\label{eq55}\dfrac{1}{2\,\lambda\,\lambda_{1}}-\dfrac{\alpha_{1}}{2\,\gamma\,\lambda_{2}}-\alpha\,\beta(1+\alpha_{1}\,\alpha_{2})=0,\end{equation}
\begin{equation}\label{eq56}-\dfrac{1}{2\,\lambda\,\lambda_{1}}+2\,\alpha(1+\alpha_{1}\,\alpha_{2})=0\end{equation}
and
\begin{equation}\label{eq57}\dfrac{\alpha_{1}}{2\,\gamma\,\lambda_{2}}-\alpha(2-\beta)(1+\alpha_{1}\,\alpha_{2})=0.\end{equation}
Moreover, from (\ref{eq48}) we deduce
 \begin{equation}\label{eq58}\dfrac{\alpha_{2}}{2\,\lambda\,\lambda_{1}}+\dfrac{1}{2\,\gamma\,\lambda_{2}}-c\,d(1+\alpha_{1}\,\alpha_{2})=0,\end{equation}
\begin{equation}\label{eq59}\dfrac{\alpha_{2}}{2\,\lambda\,\lambda_{1}}+c(2-d)(1+\alpha_{1}\,\alpha_{2})=0\end{equation}
and
\begin{equation}\label{eq60}1+\alpha_{1}\,\alpha_{2}=\dfrac{1}{4\,c\,\gamma\,\lambda_{2}}.\end{equation}
The identities (\ref{eq60}) and (\ref{eq56}) imply that $\alpha=\dfrac{c\,\gamma\,\lambda_{2}}{\lambda\,\lambda_{1}}$. From (\ref{eq59}) and (\ref{eq60}) we derive $\delta_{2}=\dfrac{(2-d)\,\lambda\,\lambda_{1}}{2\,\gamma}$. Seeing that $\delta_{1}=\alpha_{1}\,\lambda_{1}$ and $\delta_{2}=-\alpha_{2}\,\lambda_{2}$ and substituting the last identity in (\ref{eq60}) we get that $\delta_{1}=\dfrac{4\,c\,\gamma\,\lambda_{2}-1}{2\,c(2-d)\,\lambda}$.\\
Now, using (\ref{eq55}) and (\ref{eq60}) we deduce $\dfrac{\alpha\,\beta}{4\,c\,\gamma\,\lambda_{2}}=\dfrac{1}{2\,\lambda\,\lambda_{1}}-\dfrac{4\,c\,\gamma\,\lambda_{2}-1}{2\,c(2-d)\,\gamma\,\lambda\,\lambda_{1}\,\lambda_{2}},$
which reduces to $\alpha\,\beta=\dfrac{1-2\,c\,d\,\gamma\,\lambda_{2}}{(2-d)\lambda\,\lambda_{1}}$. As $\alpha=\dfrac{c\,\gamma\,\lambda_{2}}{\lambda\,\lambda_{1}}$ we derive that $\beta=\dfrac{1-2\,c\,d\,\gamma\,\lambda_{2}}{c(2-d)\gamma\,\lambda_{2}}$.\\
The solution occurs in part (5)(iii).
\\\underline{Subcase 4.4.4}: $\chi_{1}=m_{2}$ and $\chi_{2}=m_{3}$. Then, from (\ref{eq46}) we derive that $\chi_{3}=m_{1}$ because $c(2-d)\neq0$. So, using (\ref{eq48}) we get that $1+\alpha_{1}\,\alpha_{2}=-\dfrac{1}{4\,c\,\gamma\,\lambda_{2}}$. Hence, proceeding as in subcase 4.4.3 we deduce that
$$\delta_{2}=-\dfrac{(2-d)\,\lambda\,\lambda_{1}}{2\,\gamma},$$
$$\delta_{1}=-\dfrac{1+4\,c\,\gamma\,\lambda_{2}}{2\,c(2-d)\,\lambda},$$
$$\alpha=-\dfrac{c\,\gamma\,\lambda_{2}}{\lambda\,\lambda_{1}}$$
and $$\beta=\dfrac{1+4\,c\,\gamma\,\lambda_{2}}{c(2-d)\gamma\,\lambda_{2}}.$$
The solution occurs in part (5)(iv).
\\\underline{Subcase 4.4.5}: $\chi_{1}=m_{3}$ and $\chi_{2}=m_{1}$. Then, (\ref{eq46}) implies that $\chi_{3}=m_{2}$ because $c\neq0$. So that (\ref{eq47}) implies
\begin{equation}\label{eq61}\dfrac{1}{2\,\lambda\,\lambda_{1}}+2\,\alpha(1+\alpha_{1}\,\alpha_{2})=0,\end{equation}
\begin{equation}\label{eq62}-\dfrac{1}{2\,\lambda\,\lambda_{1}}-\dfrac{\alpha_{1}}{2\,\gamma\,\lambda_{2}}-\alpha\,\beta(1+\alpha_{1}\,\alpha_{2})=0\end{equation}
and
\begin{equation}\label{eq63}\dfrac{\alpha_{1}}{2\,\gamma\,\lambda_{2}}-\alpha(2-\beta)(1+\alpha_{1}\,\alpha_{2})=0.\end{equation}
Moreover, we get from (\ref{eq48})
\begin{equation}\label{eq64}\dfrac{\alpha_{2}}{2\,\lambda\,\lambda_{1}}-c\,d(1+\alpha_{1}\,\alpha_{2})=0,\end{equation}
\begin{equation}\label{eq65} -\dfrac{\alpha_{2}}{2\,\lambda\,\lambda_{1}}-c(2-d)(1+\alpha_{1}\,\alpha_{2})+\dfrac{1}{2\,\gamma\,\lambda_{2}}=0 \end{equation}
and
\begin{equation}\label{eq66}2\,c(1+\alpha_{1}\,\alpha_{2})-\dfrac{1}{2\,\gamma\,\lambda_{2}}=0.\end{equation}
From (\ref{eq66}) we get that $1+\alpha_{1}\,\alpha_{2}=\dfrac{1}{4\,c\,\gamma\,\lambda_{2}}$. Substituting this in (\ref{eq61}) we deduce $$\alpha=-\dfrac{c\,\gamma\,\lambda_{2}}{\lambda\,\lambda_{1}}.$$
Now, multiplying (\ref{eq64}) by $2$ and  (\ref{eq66}) by $d$, and adding the identities obtained we drive $\dfrac{\alpha_{2}}{\lambda\,\lambda_{1}}=\dfrac{d}{2\,\gamma\,\lambda_{2}}$. Hence, $$\delta_{2}=-\dfrac{d\,\lambda\,\lambda_{1}}{2\,\gamma}.$$
Now, since $1+\alpha_{1}\,\alpha_{2}=\dfrac{1}{4\,c\,\gamma\,\lambda_{2}}$, $\delta_{1}=\alpha_{1}\,\lambda_{1}$ and $\delta_{2}=-\alpha_{2}\,\lambda_{2}$ we derive, by using the identity above that $$\delta_{1}=\dfrac{1-4\,c\,\gamma\,\lambda_{2}}{2\,c\,d\,\lambda}.$$
On the other, using that $1+\alpha_{1}\,\alpha_{2}=\dfrac{1}{4\,c\,\gamma\,\lambda_{2}}$ and $\alpha=-\dfrac{c\,\gamma\,\lambda_{2}}{\lambda\,\lambda_{1}}$, we deduce from (\ref{eq62}) by an elementary computation  that $$\beta=\dfrac{1-2\,c(2-d)\gamma\,\lambda_{2}}{c\,d\,\gamma\,\lambda_{2}}.$$
The solution occurs in part (5)(v).
\\\underline{Subcase 4.4.6}: $\chi_{1}=m_{3}$ and $\chi_{2}=m_{2}$. Then, (\ref{eq46}) implies that $\chi_{3}=m_{1}$ because $c\neq0$. So, we derive from (\ref{eq47}) and (\ref{eq48}) the following identities
\begin{equation*}\dfrac{1}{2\,\lambda\,\lambda_{1}}+2\,\alpha(1+\alpha_{1}\,\alpha_{2})=0,\end{equation*}
\begin{equation*}-\dfrac{1}{2\,\lambda\,\lambda_{1}}+\dfrac{\alpha_{1}}{2\,\gamma\,\lambda_{2}}-\alpha(2-\beta)(1+\alpha_{1}\,\alpha_{2})=0,\end{equation*}
\begin{equation*}\dfrac{\alpha_{1}}{2\,\gamma\,\lambda_{2}}+\alpha\,\beta(1+\alpha_{1}\,\alpha_{2})=0,\end{equation*}
\begin{equation*}\dfrac{\alpha_{2}}{2\,\lambda\,\lambda_{1}}-c\,d(1+\alpha_{1}\,\alpha_{2})=0,\end{equation*}
\begin{equation*}-\dfrac{\alpha_{2}}{2\,\lambda\,\lambda_{1}}-c(2-d)(1+\alpha_{1}\,\alpha_{2})-\dfrac{1}{2\,\gamma\,\lambda_{2}}=0 \end{equation*}
and
\begin{equation*}2\,c(1+\alpha_{1}\,\alpha_{2})+\dfrac{1}{2\,\gamma\,\lambda_{2}}=0.\end{equation*}
By a similar computation to one in subcase 4.4.5 we derive that
$$\delta_{2}=\dfrac{d\,\lambda\,\lambda_{1}}{2\,\gamma},$$
$$\delta_{1}=\dfrac{1+4\,c\,\gamma\,\lambda_{2}}{2\,c\,d\,\lambda},$$
$$\alpha=\dfrac{c\,\gamma\,\lambda_{2}}{\lambda\,\lambda_{1}}$$
and
$$\beta=\dfrac{1+4\,c\,\gamma\,\lambda_{2}}{c\,d\,\gamma\,\lambda_{2}}.$$
The solution occurs in part (5)(vi).
\par Conversely if $f,g_{1},h$ and $g_{2}$ are of the
forms (1)-(5) in Theorem \ref{thm2} we check by elementary computations
that $f,g_{1},h$ and $g_{2}$ satisfy the system $(\ref{eq01})-(\ref{eq02})$, and
$f$ and $h$ are linearly independent. This completes the proof of Theorem \ref{thm2}
\end{proof}
\quad\\
\textbf{Conflict of interest} The authors declared that they have no competing interests.\\\\
\textbf{Data availability statement} Not applicable.

\end{document}